\documentclass[a4paper,11pt]{article}

\usepackage{amsmath}
\usepackage{amssymb}
\usepackage{tabularx}
\usepackage{enumerate}
\usepackage{setspace}
\usepackage{graphicx}
\usepackage{color}
\usepackage{graphicx}
\usepackage{caption}
\usepackage{subfigure}
\usepackage{bm}
\usepackage{float}

\usepackage[pagewise]{lineno}
\newtheorem{theorem}{Theorem}[section]
\newtheorem{proposition}{Proposition}[section]
\newtheorem{lemma}{Lemma}[section]
\newtheorem{remark}{Remark}[section]
\newtheorem{corollary}{Corollary}[section]
\newtheorem{defi}{Definition}[section]
\newenvironment {proof} {\noindent {\bf Proof.}}{\quad $\blacksquare$\par\vspace{3mm}}

\newcommand{\be}{\begin{equation}}
\newcommand{\ee}{\end{equation}}
\newcommand\bes{\begin{eqnarray}}
\newcommand\ees{\end{eqnarray}}
\newcommand{\bess}{\begin{eqnarray*}}
\newcommand{\eess}{\end{eqnarray*}}

\newcommand\ep{{\varepsilon}}

\begin{document}
\setlength{\baselineskip}{15.2pt} \pagestyle{myheadings}
\title{ \bf  \LARGE Eventual smoothness and sharp asymptotics in a four-dimensional Keller--Segel system with quadratic logistic degradation}
\date{\empty}
\author{
\sl{Haotian Tang $^a$, Jiashan Zheng $^b$ and Yuanyuan Ke $^a$\thanks{Corresponding author.}}\\
{ \normalsize $^a$ School of Mathematics, Renmin University of China,}\\
{ \normalsize Beijing 100872, P. R. China}\\
{ \normalsize $^b$ School of Mathematical and Informational Sciences, Yantai University,}\\
{ \normalsize Yantai, 264005, Shandong, P. R. China}\\
{}\\
}
\maketitle
\renewcommand{\thefootnote}{}
\footnotetext{Email address: haotiantang2022@163.com (H. Tang), zhengjiashan2008@163.com (J. Zheng),\\
keyy@ruc.edu.cn (Y. Ke).}
\begin{quote}\noindent {\bf Abstract.} {We study the Neumann initial-boundary problem for a fully parabolic minimal Keller-Segel system with standard logistic source in a smoothly bounded convex domain $\Omega\subset\mathbb R^4$. For arbitrary small logistic damping, we prove that the global weak solutions constructed in previous works eventually become smooth and classical, provided that the growth rate is sufficiently small. Furthermore, we derive sharp asymptotic profiles for $(u,v)$ in $L^\infty(\Omega)\times W^{1,\infty}(\Omega)$ and identify the associated leading-order coefficients in three growth regimes: for negative growth, we establish optimal exponential decay for both components, with the signal exhibiting three distinct decay patterns depending on the growth rate relative to $-1$; for zero growth, we obtain optimal algebraic decay of the solution; for suitably small positive growth, we establish the sharp exponential convergence to the positive equilibrium whenever the leading-order coefficient is nonzero, without any largeness assumption on the damping strength. These results are new, and our approach provides a valuable framework for addressing the asymptotic issues of related chemotaxis systems.}


\end{quote}

\newcommand\HI{{\bf I}}
\section{Introduction}
\indent

Chemotaxis refers to the directed movement of cells or organisms in response to chemical signals. Beyond its pivotal role in biological processes such as immune cell migration and cancer cell invasion, chemotaxis also finds broad applications in ecology, the social sciences and other fields (see \cite{P2019}).  In 1970, Keller and Segel (\cite{KS1970}) first introduced the following fully parabolic minimal chemotaxis system
\begin{equation}\label{i1111}
\left\{\aligned
&u_{t}= \Delta u-\chi\nabla\cdot(u\nabla v),&&x\in\Omega,t>0,\\
&v_{t}= \Delta v+u-v,&&x\in\Omega,t>0.
\endaligned\right.
\end{equation}
Here $u=u(x,t)$ and $v=v(x,t)$ denote the cell density and chemical signal concentration, respectively, and $\chi$ is the chemotactic sensitivity.  In one dimension, it was proved in \cite{OY2001} that the solutions to \eqref{i1111} are always globally bounded. The situation becomes considerably more intricate in higher dimensions, where various blow-up phenomena have been observed. For $N=2$, the qualitative behavior of radially symmetric solutions is closely related to the initial mass, with $\|u_{0}\|_{L^{1}(\Omega)}=\frac{8\pi}{\chi}$ representing a critical threshold. More precisely, global boundedness is guaranteed when $\|u_{0}\|_{L^{1}(\Omega)}<\frac{8\pi}{\chi}$(\cite{NSY1997}), whereas finite-time blow-up of solutions with  $\|u_{0}\|_{L^{1}(\Omega)}>\frac{8\pi}{\chi}$ has been established in \cite{HV1997,MW}. In the non-radial scenario, the corresponding critical mass decreases $\frac{4\pi}{\chi}$, although the blow-up time is still unknown (see \cite{HW2001,NSY1997,SS2001}). In dimensions $N\geq3$, the results of \cite{W2010a,W2013} further revealed that unless $\|u_{0}\|_{L^{\frac{N}{2}}(\Omega)}$ and $\|\nabla v_{0}\|_{L^{N}(\Omega)}$ are sufficiently small, finite-time blow-up of radially symmetric solutions may still occur even for arbitrarily small initial mass $\|u_{0}\|_{L^{1}(\Omega)}$.

To account for cell proliferation and competition-induced death over long time scales, we consider the chemotaxis-growth system given by
\begin{equation}\label{i22222}
\left\{\aligned
&u_{t}=\Delta u-\chi\nabla\cdot(u\nabla v)+f(u),&&x\in\Omega,t>0,\\
&v_{t}=\Delta v+u-v,&&x\in\Omega,t>0.
\endaligned\right.
\end{equation}
A typical choice is $f(u)=\lambda u-\mu u^\alpha$, where $\lambda\in\mathbb R$, $\mu>0$ and $\alpha>1$. From a mathematical perspective, the nonlinear damping term $-\mu u^\alpha$  provides additional dissipation that counteracts chemotactic aggregation. We first consider the case of quadratic damping ($\alpha=2$). In two dimensions, the basic assumption $\mu>0$ is sufficient to prevent blow-up (\cite{JX2018, L2021}). It was proved in \cite{X2018c,XZ2019} that blow-up can be ruled out even for certain sub-logistic sources, including
$f(u)=\lambda u-\frac{\mu u^2}{[\ln(u+1)]^\gamma}$ with $\gamma\in(0,1)$
and $f(u)=\lambda u-\frac{\mu u^2}{\ln\ln(u+e)}$. In dimensions $N\geq3$, the global classical solvability was established in \cite{W2010b, X2018a, X2018b, ZL2018}  for sufficiently large $\mu$. In particular, the explicit condition $\mu>\frac{N\chi}{4}$ suffices in convex
domains (\cite{W2014,X2018a}). More recently, Jiang (\cite{J2026}) showed that the same criterion remains sufficient when the second equation is replaced by $\tau v_t=\Delta v-v+u$ for any fixed $\tau>0$.  Without a largeness assumption on $\mu$, global weak solutions have been shown to exist for $N\geq3$ in \cite{L2015,ZL2018}. In a three-dimensional convex domain, Lankeit (\cite{L2015}) further established the eventual smoothness of such weak solutions for small $\lambda$, as well as their asymptotic decay when $\lambda\leq 0$. It is worth noting that whether blow-up may occur in higher dimensions ($N\geq 3$) for small positive $\mu$ remains open. Moreover, the asymptotic behavior of solutions to \eqref{i22222} has been investigated for $\lambda\in\mathbb R$ (see e.g. \cite{C2017, DW2020,W2014,X2018b}).

In the subquadratic case $1<\alpha<2$, successive studies have
extended global solvability in a generalized sense to the range $\alpha>\min\{2-\frac2N,\frac{2N+4}{N+4}\}$ (\cite{V2016,W2020b,YF2021}). Within a more general solution framework, Winkler further extended this range to all $\alpha>1$ for $N\ge2$ (\cite{W2023}) and showed that such generalized solutions converge to the positive homogeneous equilibrium in $L^1(\Omega)\times L^2(\Omega)$ when $\alpha\ge2-\frac2N$ and $\chi$ is sufficiently small
(\cite{W2020a}). Under suitable smallness conditions on $\lambda$, eventual smoothness of such generalized solutions was established in two dimensions for $\alpha>1$ (\cite{ZK2025a}), and in three dimensions first for $\alpha\in(\frac95,2)$ (\cite{VW2018}) and later for the wider range $\alpha\in(\frac53,2)$ (\cite{D2026}). The latter work also studies the asymptotic behavior of solutions in $ L^p(\Omega)\times L^\infty(\Omega)$ ($1\le p<\infty$), determining explicit convergence rates for $\lambda<0$, $\lambda=0$ and $\lambda>0$.

The asymptotic properties of global (generalized) weak solutions are of particular interest. Beyond the results recalled above (\cite{D2026, L2015, VW2018, ZK2025a}), the eventual regularity and stabilization has also been established for various chemotaxis-fluid systems in three dimensions (see e.g. \cite{BW2022, DL2023, TX2025, W2022, ZK2022}). More recently, for the basic chemotaxis--consumption system without fluid
coupling, it was demonstrated in \cite{LW2025} that the global weak solutions eventually become smooth and stabilize toward a semi-trivial steady state in general high-dimensional domains. These developments naturally lead to the following question for the minimal Keller-Segel system with signal production:

\medskip
\noindent\textbf{(Q1)}
\emph{For arbitrarily small $\mu>0$, can the smallness condition on $\lambda$ ensure the eventual smoothness of global weak solutions to \eqref{i22222} with $\alpha=2$ in $\Omega\subset\mathbb R^4$?}
\medskip

Once regularity is available, a further issue concerns the precision of the asymptotic description. The global solutions to the quadratic and subquadratic models exhibit the same qualitative large time behavior: exponential decay for $\lambda<0$, algebraic decay for $\lambda=0$, and exponential convergence to the positive constant equilibrium, provided that $\mu$ is sufficiently large. These results were established in $L^\infty(\Omega)\times L^\infty(\Omega)$ for the quadratic case (\cite{C2017,DW2020,W2014, X2018b}) and in $L^p(\Omega)\times L^\infty(\Omega)$ $1\le p<\infty$ for the subquadratic case (\cite{D2026}). At zero growth, a more precise conclusion (see \cite{CW2018,HW2019}) is available for
a Keller-Segel system with quadratic degradation and prescribed
advection: global bounded classical solutions satisfy
\begin{equation*}
\frac{c}{1+t}
\le \|u(\cdot,t)\|_{L^\infty(\Omega)}
\le \frac{C}{1+t},
\end{equation*}
which implies the the optimality of the algebraic decay rate. More recently, Zheng and Ke (\cite{ZK2025b}) also established two-sided algebraic and exponential decay estimates for a parabolic--elliptic Keller--Segel--(Navier--)Stokes system in two and three dimensions under zero and negative growth, respectively. These matching bounds, however, leave open whether the limits $\lim\limits_{t\rightarrow\infty}(1+t)\left\|u(\cdot,t)\right\|_{L^\infty(\Omega)}$  and $\lim\limits_{t\rightarrow\infty}e^{\lambda t}\left\|u(\cdot,t)\right\|_{L^\infty(\Omega)}$  exist or not.
A further question is whether sharp asymptotic profiles can
be determined for $\lambda<0$ and $\lambda>0$. This motivates a second question:

\medskip
\noindent\textbf{(Q2)}
\emph{Can one determine sharp asymptotic profiles in $L^\infty(\Omega)\times W^{1,\infty}(\Omega)$ for negative, zero and sufficiently small positive growth, identifying the leading-order coefficients?}
\medskip

The present work answers \textbf{(Q1)} and \textbf{(Q2)}
for the following Neumann initial-boundary value problem
\begin{equation}\label{model11}
\left\{\aligned
&u_{t}= \Delta u -\chi\nabla\cdot(u\nabla v)+\lambda u -\mu u ^{2},&&x\in\Omega,t>0,\\
&v_{t}= \Delta v+u-v,&&x\in\Omega,t>0,\\
&  \frac{\partial u}{\partial \nu}=\frac{\partial v}{\partial \nu}=0,&&x\in\partial\Omega,t>0,\\
&  u(x,0)=u_{0}(x),v(x,0)=v_{0}(x),&&x\in\Omega,
\endaligned\right.
\end{equation}
where $\Omega\subset\mathbb R^N$ is a domain with smooth boundary, and $\chi, \mu>0$ and $\lambda\in\mathbb R$ are fixed parameters. The unknown functions $u=u(x,t)$ and $v=v(x,t)$ represent the cell density and the chemoattractant concentration, respectively.\\
\textbf{Main results.}
We assume  that the initial data satisfy
\begin{equation}\label{initial}
\left\{\aligned
  &  u_{0}\in C^{0}(\overline{\Omega}) \ \ \mbox{with} \ \ u_{0}>0\ \ \mbox{in} \ \overline{\Omega}, \\
  &  v_{0}\in W^{1,N+1}(\Omega) \ \ \mbox{with} \ \ v_{0}\geq0\ \ \mbox{in} \ \Omega.
\endaligned\right.
\end{equation}
To facilitate our analysis, for each $\varepsilon\in(0,1)$ we consider the regularized problem
\begin{equation}\label{model111}
\left\{\aligned
&u_{\ep t}= \Delta u_{\ep} -\chi\nabla\cdot(u_{\ep}\nabla v_{\ep})+\lambda u_{\ep} -\mu u_{\ep} ^{2}-\ep u_{\ep}^\theta,&&x\in\Omega,t>0,\\
&v_{\ep t}= \Delta v_{\ep}+u_{\ep}-v_{\ep},&&x\in\Omega,t>0,\\
&  \frac{\partial u_{\ep}}{\partial \nu}=\frac{\partial v_{\ep }}{\partial \nu}=0,&&x\in\partial\Omega,t>0,\\
&  u_{\ep}(x,0)=u_{0}(x),v_{\ep}(x,0)=v_{0}(x),&&x\in\Omega
\endaligned\right.
\end{equation}
with $\theta>N+2$. The following results on global weak solvability and approximation are available for arbitrary $\mu>0$ (see \cite{L2015,ZL2018}).
\begin{proposition}\label{prop1}
Let $\Omega\subset\mathbb R^{N}$ $(N\geq3)$ be a bounded domain with smooth boundary. For any initial data fulfilling \eqref{initial} and $\mu>0$, the system \eqref{model11} admits at least one global weak solution $(u, v)$ in the sense of Definition \ref{weak solution} below. Furthermore, this solution can be chosen as a limit of solutions to the regularized problem \eqref{model111} in the following sense: for each $\varepsilon\in(0,1)$,  there exist functions
\begin{equation*}
\begin{aligned}
u_\varepsilon
&\in C^0(\overline{\Omega}\times[0,\infty))
\cap C^{2,1}(\overline{\Omega}\times(0,\infty)),\\[1mm]
v_\varepsilon
&\in C^0(\overline{\Omega}\times[0,\infty))
\cap C^{2,1}(\overline{\Omega}\times(0,\infty))\cap L^\infty_{\mathrm{loc}}
\bigl([0,\infty);W^{1,N+1}(\Omega)\bigr)
\end{aligned}
\end{equation*}
such that $u_\varepsilon>0$ and $v_\varepsilon\geq0$ in $\overline{\Omega}\times(0,\infty)$, that $(u_\varepsilon,v_\varepsilon)$ solves \eqref{model111} in the classical sense, and that for some sequence $(\varepsilon_j)_{j\in\mathbb N}\subset(0,1)$ satisfying $\varepsilon_j\searrow0$ as $j\to\infty$, one has
\begin{equation*}
u_{\varepsilon}\to u,\qquad v_{\varepsilon}\to v
\qquad\mbox{a.e. in }\Omega\times(0,\infty)
\end{equation*}
as $\varepsilon=\varepsilon_j\searrow0$.
\end{proposition}

The above existence result leaves open the question of eventual regularity of the global weak solution. Our first main result provides an affirmative answer to \textbf{(Q1)}.
\begin{theorem}\label{theorem1}
Let $\Omega\subset\mathbb R^{4}$ be a bounded convex domain with smooth boundary, and let $\mu>0$. Suppose that the initial data satisfy \eqref{initial} and $(u,v)$ denotes the global weak solution constructed in Proposition \ref{prop1}. Then there exists $\lambda_{0}>0$ such that if $\lambda<\lambda_{0}$, one can find $T_{0}>0$ such that
\begin{equation}\label{tt1}
u\in C^{2,1}(\overline{\Omega}\times[T_{0},\infty)),\ \ \ v\in C^{2,1}(\overline{\Omega}\times[T_{0},\infty)),
\end{equation}
and that the pair $(u,v)$ solves \eqref{model11} classically in $\overline{\Omega}\times[T_{0},\infty)$. Moreover, there exist $\alpha\in(0,1)$ and $C>0$ such that the solution enjoys the regularity properties
\begin{equation}\label{tt2}
\|u\|_{C^{2+\alpha,1+\frac{\alpha}{2}}(\overline{\Omega}\times[t,t+1])}+\|v\|_{C^{2+\alpha,1+\frac{\alpha}{2}}(\overline{\Omega}\times[t,t+1])}\leq C\ \ \ \mbox{for all}\ \ t\geq T_{0}.
\end{equation}
\end{theorem}

\begin{remark}\normalfont
(i) To the best of our knowledge, Theorem \ref{theorem1} provides
the first eventual smoothness result for the four-dimensional Keller--Segel system with quadratic logistic degradation.\\
(ii)
For the three-dimensional Keller--Segel system with subquadratic degradation, eventual smoothness was established for $\alpha\in(\frac95,2)$ in \cite{VW2018} and subsequently extended to $\alpha\in(\frac53,2)$ in \cite{D2026}. We expect that  a suitable adaptation of our method will further extend this result to $\alpha\in[\frac32,2)$.
\end{remark}

We next address \textbf{(Q2)} by distinguishing the three signs of $\lambda$. Throughout the asymptotic statements, $\rho_1>0$ denotes the first nonzero Neumann eigenvalue of $-\Delta$ in $\Omega$. For negative growth, we characterize the leading-order behaviour of the cell density and the signal, identifying a strictly positive coefficient for the former and distinguishing three regimes for the latter according to the position of $\lambda$ relative to $-1$.
\begin{theorem}\label{theorem2}
Let $\lambda<0$, and let
\begin{equation}\label{ttt1}
0<\delta<\min\left\{\frac{\rho_1}{12},-\lambda\right\}
\end{equation}
and
\begin{equation}\label{ttt2}
\min\{-\lambda,1\}<\sigma<\min\{\rho_1+1,\delta-\lambda\}.
\end{equation}
Under the conditions in Theorem \ref{theorem1}, there exists a constant $A>0$ such that one can find $T_{\star,1}>T_{0}$ and $C>0$ for which the following assertions hold:
\begin{enumerate}
\item[(1)]
For every $\lambda<0$, then
\begin{equation}\label{ttt3}
\|u(\cdot,t)-Ae^{\lambda t}\|_{L^\infty(\Omega)}
\leq C e^{-(\delta-\lambda)t}
\ \ \ \mbox{for all}\ t\geq T_{\star,1},
\end{equation}
and hence
\begin{equation}\label{ttt4}
\lim_{t\rightarrow\infty}\|e^{-\lambda t}u(\cdot,t)-A\|_{L^\infty(\Omega)}=0.
\end{equation}

\item[(2)] If $-1<\lambda<0$, then
\begin{equation}\label{ttt5}
\left\|v(\cdot,t)-\frac{A}{1+\lambda}e^{\lambda t}\right\|_{W^{1,\infty}(\Omega)}
\leq Ce^{-\min\{1,\sigma\}t}\ \ \ \mbox{for all}\ t\geq T_{\star,1},
\end{equation}
and hence
\begin{equation}\label{ttt6}
\lim_{t\rightarrow\infty}\left\|e^{-\lambda t}v(\cdot,t)-\frac{A}{1+\lambda}\right\|_{W^{1,\infty}(\Omega)}
=0
\end{equation}

\item[(3)] If $\lambda=-1$, then
\begin{equation}\label{ttt7}
\|v(\cdot,t)-At e^{-t}\|_{W^{1,\infty}(\Omega)}\leq Ce^{-t}
\ \ \ \mbox{for all}\ t\geq T_{\star,1},
\end{equation}
and hence
\begin{equation}\label{ttt8}
\lim_{t\rightarrow\infty}\left\|\frac{e^t}{t}v(\cdot,t)-A\right\|_{W^{1,\infty}(\Omega)}=0.
\end{equation}
\item[(4)] If $\lambda<-1$, then there exists $B>0$ such that
\begin{equation}\label{ttt9}
\|v(\cdot,t)-Be^{-t} \|_{W^{1,\infty}(\Omega)}
\leq Ce^{-\min\{\sigma,-\lambda\} t}
\ \ \ \mbox{for all}\ t\geq T_{\star,1},
\end{equation}
and hence
\begin{equation}\label{ttt10}
\lim_{t\rightarrow\infty}\|e^t v(\cdot,t)-B\|_{W^{1,\infty}(\Omega)}=0.
\end{equation}
\end{enumerate}
\end{theorem}

These asymptotic profiles further yield optimal decay estimates for the two solution components.
\begin{corollary}\label{cor1}
Under the assumptions of Theorem~\ref{theorem2}, let $A>0$ and $B>0$ be the constants obtained there. Define
\begin{equation}\label{cor12}
\Phi_\lambda(t):=
\begin{cases}
\frac{A}{1+\lambda}e^{\lambda t}, & -1<\lambda<0,\\
Ate^{-t}, & \lambda=-1,\\
Be^{-t}, & \lambda<-1.
\end{cases}
\end{equation}
For every $\epsilon_{\star}\in(0,1)$, there exists $T_{\epsilon_{\star,1}}>T_{\star,1}$ such that
\begin{equation}\label{cor11}
(1-\epsilon_{\star})Ae^{\lambda t}\leq\|u(\cdot,t)\|_{L^\infty(\Omega)}\leq(1+\epsilon_{\star})Ae^{\lambda t}\ \ \ \mbox{for all}\ t\geq T_{\epsilon_{\star,1}}
\end{equation}
and
\begin{equation}\label{cor13}
(1-\epsilon_{\star})\Phi_\lambda(t)\leq \|v(\cdot,t)\|_{W^{1,\infty}(\Omega)}
\leq(1+\epsilon_{\star})\Phi_\lambda(t)\ \ \ \mbox{for all}\ t\geq T_{\epsilon_{\star,1}},
\end{equation}
\end{corollary}
\begin{proof}
Fix $\epsilon_{\star}\in(0,1)$. By \eqref{ttt4}, \eqref{ttt6}, \eqref{ttt8} and \eqref{ttt10}, the sign-preserving property of limits along with the reverse triangle inequality provides $T_{\epsilon_{\star},1}>T_{\star,1}$ fulfilling
\begin{equation}\label{cor111}
\begin{aligned}
\left|e^{-\lambda t}\|u(\cdot,t)\|_{L^\infty(\Omega)}-A\right|
\leq&\|e^{-\lambda t}u(\cdot,t)-A\|_{L^\infty(\Omega)}\\
\leq&\epsilon_{\star}A\ \ \ \mbox{for all}\ t\geq T_{\epsilon_{\star,1}}
\end{aligned}
\end{equation}
and
\begin{equation}\label{cor112}
\begin{aligned}
\left|\frac{\|v(\cdot,t)\|_{W^{1,\infty}(\Omega)}}{\Phi_\lambda(t)}-1\right|
\leq&\left\|\frac{v(\cdot,t)}{\Phi_\lambda(t)}-1\right\|_{W^{1,\infty}(\Omega)}\\
\leq&\epsilon_{\star}\ \ \ \mbox{for all}\ t\geq T_{\epsilon_{\star,1}},
\end{aligned}
\end{equation}
which are equivalent to the inequality in \eqref{cor11} and \eqref{cor13}.
\end{proof}

\begin{remark}\normalfont
(i) To our knowledge, these are the first optimal exponential decay estimates for the four-dimensional Keller--Segel system \eqref{model11} under negative growth.\\
(ii) Theorem~\ref{theorem2} establishes asymptotic expansions
for both solution components with strictly positive
leading-order coefficients and provides explicit estimates for the corresponding remainders in
$L^\infty(\Omega)\times W^{1,\infty}(\Omega)$.
In particular, Corollary~\ref{cor1} yields
\begin{equation}
    \lim_{t\to\infty}
    e^{-\lambda t}\|u(\cdot,t)\|_{L^\infty(\Omega)}
    = A, \qquad \lambda<0
\end{equation}
and as $t\to\infty$,
\begin{equation}
    \|v(\cdot,t)\|_{W^{1,\infty}(\Omega)}
    \sim
    \begin{cases}
        \displaystyle\frac{A}{1+\lambda}e^{\lambda t},
        & -1<\lambda<0,\\[6pt]
        Ate^{-t}, & \lambda=-1,\\[3pt]
        Be^{-t}, & \lambda<-1.
    \end{cases}
\end{equation}
(iii) While the cell density has the asymptotic decay scale
$e^{\lambda t}$ for all $\lambda<0$, the signal
concentration exhibits three distinct asymptotic
decay regimes: $e^{\lambda t}$ for $-1<\lambda<0$,
$te^{-t}$ for $\lambda=-1$, and $e^{-t}$ for
$\lambda<-1$. The threshold $-\lambda=1$ reflects reflects the competition
between the cell decay rate $-\lambda$ and the intrinsic
signal relaxation rate $1$.
\end{remark}

In the critical case $\lambda=0$, we establish the exponential convergence of $u$ towards an explicit algebraic profile and derive a second-order asymptotic expansion for $v$.
\begin{theorem}\label{theorem3}
Let $\lambda=0$. Under the assumptions in Theorem \ref{theorem1}, there exist $C_\infty\in\mathbb R$, $T_{\star,2}>T_{0}$ and $C>0$ such that
\begin{equation}\label{tttt1}
\left\|u(\cdot,t)-\frac1{\mu t+C_\infty}\right\|_{L^\infty(\Omega)}
\leq Ce^{-\frac{\rho_1}{12}t}\ \ \ \mbox{for all}\ t\geq T_{\star,2}
\end{equation}
and
\begin{equation}\label{tttt2}
\left\|v(\cdot,t)-\frac1{\mu t+C_\infty}
-\frac\mu{(\mu t+C_\infty)^2}\right\|_{W^{1,\infty}(\Omega)}
\leq\frac C{(1+t)^3}\ \ \ \mbox{for all}\ t\geq T_{\star,2}.
\end{equation}
Furthermore, one has
\begin{equation}\label{tttt3}
\lim_{t\rightarrow\infty}\left\|(1+t)u(\cdot,t)-\frac1{\mu }\right\|_{L^\infty(\Omega)}=0
\end{equation}
and
\begin{equation}\label{tttt4}
\lim_{t\rightarrow\infty}\left\|(1+t)v(\cdot,t)-\frac1{\mu }\right\|_{W^{1,\infty}(\Omega)}=0.
\end{equation}
\end{theorem}
In particular, we obtain optimal algebraic decay rates for u and v.
\begin{corollary}\label{cor2}
Suppose that the conditions of Theorem \ref{theorem3} hold. Then for every $\epsilon_{\star}\in(0,1)$, there exists $T_{\epsilon_{\star},2}>T_{\star,2}$ such that
\begin{equation}\label{cor21}
\frac{1-\epsilon_{\star}}{\mu(1+t)}\leq\|u(\cdot,t)\|_{L^\infty(\Omega)}\leq\frac{1+\epsilon_{\star}}{\mu(1+t)}\ \ \ \mbox{for all}\ t\geq T_{\epsilon_{\star,2}}
\end{equation}
and
\begin{equation}\label{cor22}
\frac{1-\epsilon_{\star}}{\mu(1+t)}\leq\|v(\cdot,t)\|_{W^{1,\infty}(\Omega)}\leq\frac{1+\epsilon_{\star}}{\mu(1+t)}\ \ \ \mbox{for all}\ t\geq T_{\epsilon_{\star,2}}.
\end{equation}
\end{corollary}

\begin{proof}
Fix $\epsilon_{\star}\in(0,1)$. According to \eqref{tttt3} and \eqref{tttt4}, we use the reverse triangle inequality and the sign-preserving property of limits to find $T_{\epsilon_{\star},2}>T_{\star,2}$ such that
\begin{equation}\label{cor221}
\begin{aligned}
\left|(1+t)\|u(\cdot,t)\|_{L^\infty(\Omega)}-\frac1\mu\right|
\leq&\left\|(1+t)u(\cdot,t)-\frac1\mu\right\|_{L^\infty(\Omega)}\\
\leq&\frac{\epsilon_{\star}}\mu\ \ \mbox{for all}\ t\geq T_{\epsilon_{\star,2}}
\end{aligned}
\end{equation}
and
\begin{equation}\label{cor222}
\begin{aligned}
\left|(1+t)\|v(\cdot,t)\|_{W^{1,\infty}(\Omega)}-\frac1\mu\right|
\leq&\left\|(1+t)v(\cdot,t)-\frac1\mu\right\|_{W^{1,\infty}(\Omega)}\\
\leq&\frac{\epsilon_{\star}}\mu\ \ \mbox{for all}\ t\geq T_{\epsilon_{\star,2}},
\end{aligned}
\end{equation}
and thus establishes \eqref{cor21} and \eqref{cor22}.
\end{proof}

\begin{remark}\normalfont
(i) The two-sided algebraic estimates in \cite{CW2018,HW2019,ZK2025b} imply that
\begin{equation*}
    \lim_{t\to\infty}
    (1+t)^\ell\|u(\cdot,t)\|_{L^\infty(\Omega)}=0
    \qquad\text{for every }\ell\in(0,1).
\end{equation*}
At the critical exponent $\ell=1$, however, such estimates only provide positive upper and lower bounds and do not ensure
convergence. Theorem~\ref{theorem3} takes a further step by establishing the following precise limits in terms of the strong norms of $L^\infty(\Omega)$ and $W^{1,\infty}(\Omega)$:
\begin{equation*}
    \lim_{t\to\infty}
    (1+t)\|u(\cdot,t)\|_{L^\infty(\Omega)}
    =\frac{1}{\mu}
    \qquad\text{and}\qquad
    \lim_{t\to\infty}
    (1+t)\|v(\cdot,t)\|_{W^{1,\infty}(\Omega)}
    =\frac{1}{\mu}.
\end{equation*}
(ii) The question of whether chemotactic cross-diffusion influences the fine structure of large time asymptotics was  raised after Theorem 1.1 in \cite{CW2018}.
In the present fluid-free setting, Theorem~\ref{theorem3}
identifies an explicit leading coefficient $\frac{1}{\mu}$
independent of the chemotactic sensitivity and the initial
data, together with a solution-dependent constant $C_\infty$.\\
(iii) In the absence of proliferation, the finer asymptotic structure is characterized by the common leading term $(\mu t+C_\infty)^{-1}$ for both components and an additional second-order correction $\mu(\mu t+C_\infty)^{-2}$ in the asymptotic expansion of $v$.
\end{remark}

For sufficiently small $\lambda>0$, we show refined asymptotic expansions of $u$ and $v$ about the spatially homogeneous equilibrium $\left(\frac{\lambda}{\mu},\frac{\lambda}{\mu}\right)$. In contrast to the case $\lambda<0$, the corresponding asymptotic coefficient $\widetilde{A}$ may vanish. Notably, no additional restriction on $\mu$ is imposed.
\begin{theorem}\label{theorem4}
Let the conditions in Theorem \ref{theorem1} hold. Then there exists $\tilde{\lambda}_{0}>0$ with the property that if $0<\lambda<\tilde{\lambda}_{0}$, one can find $\widetilde{A}\in\mathbb R$, $T_{\star,3}>T_{0}$ and $C>0$ such that
\begin{equation}\label{ttttt1}
\left\|u(\cdot,t)-\frac\lambda\mu-\widetilde{A}e^{-\lambda t}
\right\|_{L^\infty(\Omega)}\leq Ce^{-2\lambda t}
\ \ \ \mbox{for all}\ t\geq T_{\star,3}
\end{equation}
and
\begin{equation}\label{ttttt2}
\left\|v(\cdot,t)-\frac\lambda\mu-\frac{\widetilde{A}}{1-\lambda}e^{-\lambda t}\right\|_{W^{1,\infty}(\Omega)}\leq Ce^{-2\lambda t}
\ \ \ \mbox{for all}\ t\geq T_{\star,3}.
\end{equation}
Furthermore, we have
\begin{equation}\label{ttttt3}
\lim_{t\rightarrow\infty}\left\|e^{\lambda t}\left(u(\cdot,t)-\frac\lambda\mu\right)-\widetilde{A}\right\|_{L^\infty(\Omega)}=0
\end{equation}
and
\begin{equation}\label{ttttt4}
\lim_{t\rightarrow\infty}\left\|e^{\lambda t}\left(v(\cdot,t)-\frac\lambda\mu\right)-\frac{\widetilde{A}}{1-\lambda}\right\|_{W^{1,\infty}(\Omega)}=0.
\end{equation}
\end{theorem}

The exponential convergence rate is optimal whenever $\widetilde{A}\neq0$, as shown in the following corollary.
\begin{corollary}\label{cor3}
Under the hypotheses of Theorem~\ref{theorem4}, let $\widetilde{A}\in\mathbb R$ and $\tilde{\lambda}_{0}$ be the constants obtained there. Suppose that $0<\lambda<\tilde{\lambda}_{0}$. If $\widetilde{A}\neq0$, then for every $\epsilon_{\star}\in(0,1)$ there exists $T_{\epsilon_{\star},3}>T_{\star,3}$ such that
\begin{equation}\label{cor31}
(1-\epsilon_{\star})|\widetilde{A}|e^{-\lambda t}\leq\left\|u(\cdot,t)-\frac\lambda\mu\right\|_{L^\infty(\Omega)}\leq(1+\epsilon_{\star})|\widetilde{A}|e^{-\lambda t}\ \ \mbox{for all}\ t\geq T_{\epsilon_{\star,3}}
\end{equation}
and
\begin{equation}\label{cor32}
\frac{(1-\epsilon_{\star})|\widetilde{A}|}{1-\lambda }e^{-\lambda t}\leq\left\|v(\cdot,t)-\frac\lambda\mu\right\|_{W^{1,\infty}(\Omega)}\leq\frac{(1+\epsilon_{\star})|\widetilde{A}|}{1-\lambda }e^{-\lambda t}\ \ \mbox{for all}\ t\geq T_{\epsilon_{\star,3}}.
\end{equation}
\end{corollary}

\begin{proof}
Let $\widetilde A\neq0$ and $\epsilon_{\star}\in(0,1)$.
In view of \eqref{ttttt3} and \eqref{ttttt4}, the sign-preserving property of limits together with the reverse triangle inequality yields the existence of $T_{\epsilon_{\star},3}>T_{\star,3}$ such that
\begin{equation}\label{cor331}
\begin{aligned}
\left|e^{\lambda t}\left\|u(\cdot,t)-\frac\lambda\mu\right\|_{L^\infty(\Omega)}
-|\widetilde A|\right|\leq&\left\|e^{\lambda t}\left(u(\cdot,t)-\frac\lambda\mu\right)-\widetilde{A}\right\|_{L^\infty(\Omega)}\\
\leq&\epsilon_{\star}|\widetilde A|\ \ \mbox{for all}\ t\geq T_{\epsilon_{\star,3}}
\end{aligned}\end{equation}
and
\begin{equation}\label{cor332}
\begin{aligned}
\left|e^{\lambda t}\left\|v(\cdot,t)-\frac\lambda\mu\right\|_{W^{1,\infty}(\Omega)}
-\frac{|\widetilde A|}{1-\lambda}\right|\leq&\left\|e^{\lambda t}\left(v(\cdot,t)-\frac\lambda\mu\right)-\frac{\widetilde{A}}{1-\lambda}\right\|_{W^{1,\infty}(\Omega)}\\
\leq&\frac{\epsilon_{\star}|\widetilde A|}{1-\lambda}\ \ \mbox{for all}\ t\geq T_{\epsilon_{\star,3}}.
\end{aligned}
\end{equation}
This readily proves \eqref{cor31} and \eqref{cor32}.
\end{proof}

\begin{remark}\normalfont
(i) For positive growth, earlier studies establish stabilization or exponential upper bounds under their respective damping and boundedness assumptions
(see e.g. \cite{C2017,D2026,DW2020,DX2022,HZ2016,W2014,X2018b,ZK2025b}). Theorem~\ref{theorem4} further provides precise
asymptotic expansions around the positive equilibrium in the strong norms of $L^\infty(\Omega)$ and
$W^{1,\infty}(\Omega)$, identifying the respective coefficients of $e^{-\lambda t}$ as $\widetilde A$
and $\frac{\widetilde A}{1-\lambda}$.\\
(ii)When $\widetilde A\neq0$, Corollary~\ref{cor3} yields the optimality of the convergence rate $e^{-\lambda t}$ for both components. To the best of our knowledge, this is the first
sharp asymptotic description for the Keller-Segel system \eqref{model11} without any largeness assumption on $\mu$,
provided that $\lambda>0$ is sufficiently small.\\
(iii) A sufficient condition for $\widetilde A\neq0$ is
\begin{equation*}
\frac{1}{|\Omega|}\int_\Omega u_0<\frac{\lambda}{\mu},
\end{equation*}
meaning that the initial mean cell density lies below its positive equilibrium value.
\end{remark}

\begin{remark}\normalfont
As noted by Bellomo et al.~\cite[p.~1691]{BBTW2015},
results on the convergence of solutions to a single equilibrium were largely restricted to the case
of sufficiently strong logistic damping. In Remark~(iv) of \cite{W2022}, Winkler also raises the question of whether the large time behaviour of solutions can be described more precisely beyond eventual boundedness.
Theorems~\ref{theorem2}--\ref{theorem4} address this question for the present four-dimensional fluid-free system \eqref{model11}. More precisely, for arbitrary $\mu>0$, we obtain refined asymptotic profiles of both components in $L^\infty(\Omega)\times W^{1,\infty}(\Omega)$ for $\lambda<0$, $\lambda=0$, and sufficiently small $\lambda>0$.
\end{remark}

\begin{remark}\normalfont
For spatially homogeneous initial data, \eqref{model11} reduces to the ODE system
\begin{equation*}
\begin{cases}
U'=\lambda U-\mu U^2,\\
V'+V=U,\\
U(0)=U_0>0,\quad V(0)=V_0\geq0.
\end{cases}
\end{equation*}
Its solution is explicitly given by
\begin{equation*}
\left\{
\begin{aligned}
U(t)&=
\begin{cases}
\displaystyle\frac{1}{\mu t+U_0^{-1}},
& \lambda=0,\\[8pt]
\displaystyle\frac{1}{\frac{\mu}{\lambda}+\left(U_0^{-1}-\frac{\mu}{\lambda}\right)
e^{-\lambda t}},
& \lambda\neq0,
\end{cases}\\[8pt]
V(t)&=e^{-t}V_0+\int_0^t e^{-(t-s)}U(s)\,ds.
\end{aligned}
\right.
\end{equation*}
These formulas provide a natural benchmark for the asymptotic profiles established in Theorems~\ref{theorem2}--\ref{theorem4}.

For $\lambda<0$, as $t\to\infty$, we have
\begin{equation*}
\left\{
\begin{aligned}
U(t)&\sim ae^{\lambda t},\\[4pt]
V(t)&\sim
\begin{cases}
\displaystyle\frac{a}{1+\lambda}e^{\lambda t},
& -1<\lambda<0,\\[6pt]
ate^{-t}, & \lambda=-1,\\[2pt]
be^{-t}, & \lambda<-1,
\end{cases}
\end{aligned}
\right.
\end{equation*}
where
$a=(U_0^{-1}-\mu/\lambda)^{-1}>0$ and in the case $\lambda<-1$, $b=V_0+\int_0^\infty e^sU(s)\,ds>0$.
Thus, the three distinct decay regimes for the signal in Theorem~\ref{theorem2}, including the additional factor $t$ at $\lambda=-1$, are already present in the spatially homogeneous dynamics.

For $\lambda=0$, we have
\begin{equation*}
\left\{
\begin{aligned}
U(t)
&=\frac{1}{\mu t+U_0^{-1}},\\[6pt]
V(t)
&=\frac{1}{\mu t+U_0^{-1}}
+\frac{\mu}{(\mu t+U_0^{-1})^2}
+O(t^{-3}).
\end{aligned}
\right.
\end{equation*}
Theorem~\ref{theorem3} recovers precisely this asymptotic structure, with $U_0^{-1}$ replaced by the effective shift $C_\infty$. In particular, both components retain the same universal leading coefficient $\frac{1}{\mu}$.

For $0<\lambda<\frac12$, we have
\begin{equation*}
\left\{
\begin{aligned}
U(t)
&=\frac{\lambda}{\mu}
+\widetilde a e^{-\lambda t}
+O(e^{-2\lambda t}),\\[6pt]
V(t)
&=\frac{\lambda}{\mu}
+\frac{\widetilde a}{1-\lambda}e^{-\lambda t}
+O(e^{-2\lambda t}),
\end{aligned}
\right.
\end{equation*}
where $\widetilde a=\frac{\lambda}{\mu}
\left(1-\frac{\lambda}{\mu U_0}\right)$. Theorem~\ref{theorem4} yields the same rate of convergence and the same relation between the
leading coefficients of $u$ and $v$. Moreover, the rate $e^{-\lambda t}$ is optimal whenever the leading coefficient is nonzero.

Taken together, these results demonstrate that the asymptotic structure of the spatially homogeneous dynamics persists for the spatially inhomogeneous
weak solutions considered here in $L^\infty(\Omega)\times W^{1,\infty}(\Omega)$. However, the constants $A$, $B$, $C_\infty$
and $\widetilde A$ depend on the spatial evolution of the solutions and need not coincide with those of the homogeneous system with initial data given by the spatial averages of $u_0$ and $v_0$.
\end{remark}
\textbf{Mathematical challenges.}
Two difficulties arise in addressing \textbf{(Q1)}.
The first is to establish eventual smallness of the coupled quantity
\begin{equation}\label{quantity}
\int_\Omega u_\varepsilon^2
+\int_\Omega|\nabla v_\varepsilon|^4.
\end{equation}
The three-dimensional energy argument in \cite{L2015} relies
on controlling $\int_\Omega|\nabla v_\varepsilon|^6$ by an
arbitrarily small multiple of
$\int_\Omega|\nabla|\nabla v_\varepsilon|^2|^2$,
together with powers of $\int_\Omega|\nabla v_\varepsilon|^4$.
In four dimensions, however, the interpolation becomes critical, so that Young's inequality no longer yields the desired absorption without an additional smallness assumption. Starting from arbitrary initial data,
one must establish that the energy becomes sufficiently
small and remains so, with estimates and a waiting time
independent of $\varepsilon$. The second difficulty concerns the passage from this eventual energy control to higher regularity. In three dimensions, the $L^2(\Omega)\times L^4(\Omega)$ bounds for $(u_\varepsilon,\nabla v_\varepsilon)$ ensure higher integrability than in the critical space
$L^{\frac{N}{2}}(\Omega)\times L^N(\Omega)$. By contrast, both bounds are exactly critical for $N=4$ (\cite{C2015,J2020}). At this critical level, the direct semigroup argument used in \cite{L2015} encounters a nonintegrable time singularity. An additional gain in integrability is thus needed to obtain eventual boundedness in $L^\infty(\Omega)\times W^{1,\infty}(\Omega)$ and subsequently establish classical regularity.

For \textbf{(Q2)}, the challenge is to determine precise leading-order asymptotic profiles beyond the information supplied by stabilization and decay estimates. For $\lambda>0$, exponential convergence to the positive equilibrium is known under sufficiently strong logistic damping (\cite{C2017,DW2020,W2014}). For $\lambda=0$, matching algebraic upper and lower bounds have been obtained under quadratic degradation, even in the presence of a prescribed fluid flow (\cite{CW2018,HW2019}). For $\lambda<0$, more recent results establish two-sided exponential decay bounds for parabolic--elliptic Keller--Segel--(Navier--)Stokes systems in two and three
dimensions (\cite{ZK2025b}). These results, however, do not provide the precise asymptotics sought here for \eqref{model11} in four dimensions for any fixed $\mu>0$. Indeed, even matching decay bounds do not ensure convergence after rescaling or identify the limiting coefficients, while an
exponential upper estimate alone does not determine the exact decay rate. A further difficulty lies in deriving the exact asymptotic behaviors in the topology of $L^\infty(\Omega)\times W^{1,\infty}(\Omega)$.\\
\textbf{Main ideas and plan of the paper.}
In Section 2, we collect the preliminary definitions and several useful estimates. In Section 3, we retain the coupled energy \eqref{quantity} used in \cite{L2015}. After deriving the basic estimates in Lemma \ref{d}, we combine H\"older's inequality with
the four-dimensional Sobolev embedding
$W^{1,2}(\Omega)\hookrightarrow L^4(\Omega)$ to establish the interpolation inequality
\begin{equation*}
\int_\Omega|\nabla v_\varepsilon|^6
\leq C_\star
\left(\int_\Omega|\nabla v_\varepsilon|^4\right)^{\frac12}
\left(
\int_\Omega\left|\nabla|\nabla v_\varepsilon|^2\right|^2
+\int_\Omega|\nabla v_\varepsilon|^4
\right)
\end{equation*}
(see Lemma \ref{a}). This estimate enables us to absorb the nonlinear gradient term into the dissipation whenever the coupled quantity is sufficiently small. Building on this observation, we employ a first self-map type argument to show that the smallness of \eqref{quantity} persists once attained, provided $\lambda_+$ is sufficiently small (Lemmas  \ref{b}-\ref{c}). Subsequently, we show the eventual smallness of \eqref{quantity} after some waiting time independent of $\varepsilon$ in Proposition \ref{e}.

With this uniform smallness available, a second self-mapping argument raises the integrability beyond the critical exponents. For suitably large
time $t_0$, we consider
\begin{equation*}
(t-t_0)^{\frac13}\|u_\varepsilon(\cdot,t)\|_{L^3(\Omega)}
\quad\text{and}\quad
(t-t_0)^{\frac16}\|\nabla v_\varepsilon(\cdot,t)\|_{L^6(\Omega)}.
\end{equation*}
The time weights and exponents are chosen to control the singularities arising in the semigroup estimates. Relying on the smoothing $L^{p}$--$L^{q}$ estimates for the Neumann heat semigroup, we can close the loop argument and obtain the eventual smallness of $\|u_\varepsilon(\cdot,t)\|_{L^3(\Omega)}$ and $\|\nabla v_\varepsilon(\cdot,t)\|_{L^6(\Omega)}$ in Lemma \ref{f}. Further semigroup estimates yield eventual uniform bounds in
$L^\infty(\Omega)\times W^{1,\infty}(\Omega)$ (see Lemma \ref{g}).
These bounds together with standard parabolic regularity and
compactness theory ensure the eventual regularity of the global weak solutions (Lemmas \ref{i}-\ref{ii}). Moreover, by suitably reducing the growth threshold, the eventual
$L^\infty(\Omega)$ bound for $u$ can be made as small as needed (see Corollary \ref{iii}), which is crucial for the subsequent asymptotic analysis.

In Section 4, our approach combines estimates for the deviations of the
solutions from their spatial averages with an analysis of the evolution of these averages. The key preliminary results are established in Lemmas~\ref{j}--\ref{m}.
We first exclude loss of positive cell mass in the approximation
limit and establish positivity in Lemma~\ref{j}, then derive
the exact evolution equations for $\bar u$ and $\bar v$ in Lemma~\ref{k}.
The arbitrarily small eventual $L^\infty(\Omega)$ bound for
$u$ enables us to derive coupled energy estimates, which along with interpolation and semigroup estimates yield exponential decay of $\|u-\bar u\|_{L^\infty(\Omega)}$ and $\|v-\bar v\|_{W^{1,\infty}(\Omega)}$ (see Lemma~\ref{l}). For zero and sufficiently small positive growth,
we further exploit a lower bound for $\bar u$ in conjunction
with the preceding energy estimates to establish exponential
decay of the relative variance
\begin{equation*}
\frac{\frac1{|\Omega|}\int_\Omega(u-\bar u)^2}{\bar u(t)^2}
\end{equation*}
in Lemma~\ref{m}.
This controls the error term in the equation for $\frac{1}{\bar u}$
and prepares the subsequent analysis of the spatial averages.

In Section 4.1, we consider the case $\lambda<0$.
Interpolation first gives the time integrability of
$\|\Delta v(\cdot,t)\|_{L^\infty(\Omega)}$.
A subsequent comparison argument gives
$\|u(\cdot,t)\|_{L^\infty(\Omega)}=O(e^{\lambda t})$,
which enables us to establish the convergence of
$e^{-\lambda t}\bar u(t)$ to a strictly positive constant.
A further semigroup estimate for $e^{-\lambda t}(u-\bar u)$
then yields the asymptotic profile of $u$. The corresponding profile of $v$ is derived from its averaged equation and spatial fluctuation estimates, with separate treatments of the three regimes $-1<\lambda<0$, $\lambda=-1$
and $\lambda<-1$, thereby completing the proof of Theorem~\ref{theorem2}.

The arguments in Sections 4.2 and 4.3 both rely on the equation
for the reciprocal of the spatial mean,
\begin{equation*}
\frac{d}{dt}\frac{1}{\bar u(t)}
=-\frac{\lambda}{\bar u(t)}+\mu+\mu\frac{\frac1{|\Omega|}\int_\Omega(u-\bar u)^2}{\bar u^2(t)},
\end{equation*}
whose last term is controlled by Lemma~\ref{m}. When $\lambda=0$, this estimate ensures that the derivative of $\frac{1}{\bar u(t)}-\mu t$ is integrable, so that this quantity converges to a finite limit $C_\infty$.
Together with the decay of the spatial deviations, this proves the exponential convergence of $u$ towards the asymptotic profile $(\mu t+C_\infty)^{-1}$. Substituting this profile into the averaged signal equation, we identify the correction term $\mu(\mu t+C_\infty)^{-2}$ and obtain the remainder estimate in Theorem~\ref{theorem3}. For small $\lambda>0$, solving the linear equation for
$\frac{1}{\bar u}-\frac{\mu}{\lambda}$ determines the
coefficient $\widetilde A$ of the term $e^{-\lambda t}$
in the expansion of $\bar u$. Using the identity
\begin{equation*}
\bar u(t)-\frac{\lambda}{\mu}-\widetilde A e^{-\lambda t}
=-\frac{\lambda^2}{\mu^2}\left[\frac{1}{\bar u(t)}-\frac{\mu}{\lambda}
+\frac{\mu^2\widetilde A}{\lambda^2}e^{-\lambda t}\right]
+\frac{\lambda^2}{\mu^2}\bar u(t)\left[\frac{1}{\bar u(t)}-\frac{\mu}{\lambda}\right]^2,
\end{equation*}
we arrive at the desired expansion for $\bar u$ with an
$O(e^{-2\lambda t})$ remainder.
In conjunction with the faster decay of the spatial deviations,
this yields the asserted expansion for $u$. Once the asymptotic profile of $u$ is established, the averaged signal equation and the decay of $v-\bar v$ yield the corresponding expansion for $v$, with $\frac{\widetilde A}{1-\lambda}$ as the coefficient of $e^{-\lambda t}$. This completes the proof of Theorem~\ref{theorem4}.\\
\textbf{Notations.} Throughout this paper, $C$, $C_{\star}$ and $C_{i}$ $(i=1, 2, \cdot\cdot\cdot)$ denote various positive constants
independent of time. We set $n_{+}:=\max\{0,n\}$ and  $n_{-}:=\max\{-n,0\}$. For simplicity, we omit ${\rm d}x$ in all spatial integrals.

\section{Preliminaries}
\setcounter{equation}{0}
\indent

In this section, we recall the definition of a global weak solution and some auxiliary results needed in the proofs below. We begin with the following definition.
\begin{defi}\label{weak solution}
Let $u_0$ and $v_0$ satisfy \eqref{initial}. A pair $(u,v)$ of nonnegative functions
\begin{equation}\label{dd1}
\left\{\begin{aligned}
u&\in L^2_{\mathrm{loc}}([0,\infty);L^2(\Omega)),\\
v&\in L^1_{\mathrm{loc}}([0,\infty);W^{1,1}(\Omega)),
\end{aligned}\right.
\end{equation}
is called a global weak solution of \eqref{model11} if
\begin{equation}\label{dd2}
\nabla u\ \ \mbox{and}\ \ u \nabla v \ \ \mbox{belong to} \ \ L^{1}_{loc}(\overline{\Omega}\times[0,\infty);\mathbb{R}^{N}),
\end{equation}
and if the identities
\begin{equation}\label{dd3}
\begin{aligned}
-\int_0^\infty\!\int_\Omega u\varphi_t
-\int_\Omega u_0\varphi(\cdot,0)
={}&-\int_0^\infty\!\int_\Omega\nabla u\cdot\nabla\varphi
+\chi\int_0^\infty\!\int_\Omega u\nabla v\cdot\nabla\varphi+\int_0^\infty\!\int_\Omega(\lambda u-\mu u^2)\varphi
\end{aligned}
\end{equation}
and
\begin{equation}\label{dd4}
\begin{aligned}
-\int_0^\infty\!\int_\Omega v\varphi_t
-\int_\Omega v_0\varphi(\cdot,0)
={}&-\int_0^\infty\!\int_\Omega\nabla v\cdot\nabla\varphi-\int_0^\infty\!\int_\Omega v\varphi
+\int_0^\infty\!\int_\Omega u\varphi
\end{aligned}
\end{equation}
hold for all $\varphi\in C_c^\infty(\overline{\Omega}\times[0,\infty))$.
\end{defi}

We next recall the sharp Poincar\'e inequality, which follows
directly from the variational characterization of the first
nonzero Neumann eigenvalue, as detailed in \cite{WX2012}.
\begin{lemma}\label{p1}
Let $\Omega\subset\mathbb R^{N}$ be a bounded domain with smooth boundary, and let $\rho_1$ denote the first nonzero eigenvalue of $-\Delta$ in $\Omega$ subject to homogeneous Neumann boundary conditions. Define $\overline{w}=\frac{1}{|\Omega|}\int_{\Omega}w$. Then
\begin{equation}\label{poincare1}
\int_{\Omega}(w-\overline{w})^{2}\leq\frac{1}{\rho_1}\int_{\Omega}|\nabla w|^{2}\ \ \ \mbox{for all}\ w\in W^{1,2}(\Omega),
\end{equation}
Equivalently,
\begin{equation}\label{poincare2}
\int_{\Omega}w^{2}\leq\frac{1}{\rho_1}\int_{\Omega}|\nabla w|^{2}+|\Omega|\overline{w}^{2}\ \ \ \mbox{for all}\ w\in W^{1,2}(\Omega).
\end{equation}
The constant $\rho_1^{-1}$ is optimal. In what follows, we may also write $C_P$ for $\rho_1^{-1}$.
\end{lemma}

We shall also use the following standard estimates for the Neumann heat semigroup. For their proofs, we refer to \cite{C2015,CL2016,CLX2025,W2010a}.
\begin{lemma}\label{p2}
Let $\{e^{t\Delta}\}_{t\geq0}$ denote the Neumann heat semigroup in $\Omega$, and let $\rho_1$ denote the first nonzero eigenvalue of $-\Delta$ in $\Omega$ under Neumann boundary conditions. Then there exist positive constants $k_i$ $(i=1, ..., 4)$ and $\widehat{k}_1$, depending only on $\Omega$ such that
\begin{enumerate}
\item[(1)] If $1\leq p\leq q\leq\infty$, then
\begin{equation}\label{p2a}
\|e^{t\Delta}w\|_{L^q(\Omega)}\leq k_1\left(1+t^{-\frac N2(\frac1p-\frac1q)}\right)e^{-\rho_1t}\|w\|_{L^p(\Omega)}, \ \ \forall t>0
\end{equation}
holds for all $w\in L^p(\Omega)$ satisfying $\int_\Omega w=0$.
\item[(2)] If $1\leq p\leq q\leq\infty$, then
\begin{equation}\label{p2e}
\|e^{t\Delta}w\|_{L^q(\Omega)}\leq\widehat{k}_1\left(1+t^{-\frac N2(\frac1p-\frac1q)}\right)\|w\|_{L^p(\Omega)}, \ \ \forall t>0
\end{equation}
holds for all $w\in L^p(\Omega)$.
\item[(3)] If $1\leq p\leq q\leq\infty$, then
\begin{equation}\label{p2b}
\|\nabla e^{t\Delta}w\|_{L^q(\Omega)}\leq k_2\left(1+t^{-\frac12-\frac N2(\frac1p-\frac1q)}\right)e^{-\rho_1t}\|w\|_{L^p(\Omega)}, \ \ \forall t>0
\end{equation}
is valid for all $w\in L^p(\Omega)$.
\item[(4)] If $2\leq p\leq q\leq\infty$, then
\begin{equation}\label{p2c}
\|\nabla e^{t\Delta}w\|_{L^q(\Omega)}\leq k_3\left(1+t^{-\frac N2(\frac1p-\frac1q)}\right)e^{-\rho_1t}\|\nabla w\|_{L^p(\Omega)}, \ \ \forall t>0
\end{equation}
is true for all $w\in W^{1,q}(\Omega)$.
\item[(5)] If $1<p\leq q\leq\infty$, then
\begin{equation}\label{p2d}
\|e^{t\Delta}\nabla\cdot w\|_{L^q(\Omega)}\leq k_4\left(1+t^{-\frac12-\frac N2(\frac1p-\frac1q)}\right)e^{-\rho_1t}\|w\|_{L^p(\Omega)}, \ \ \forall t>0
\end{equation}
holds for all $w\in(L^p(\Omega))^N$.
\end{enumerate}
\end{lemma}

\section{Eventual smoothness in 4D}
\setcounter{equation}{0}
\indent

This section is devoted to the proof of Theorem~1.1.
Our aim is to establish the eventual smoothness of the global
weak solution in four dimensions, provided that the growth
rate is sufficiently small.
We begin with some basic $\varepsilon$-independent estimates
for the approximate solutions, which will be used in the
subsequent analysis.
\begin{lemma}\label{d}
Let $\hat\lambda>\lambda_+$. Then for each $\varepsilon\in(0,1)$, the solution to \eqref{model111} has the following properties:
\begin{enumerate}
\item[(1)]
\begin{equation}\label{d0}
\int_{\Omega}u_{\varepsilon}(\cdot,t)\leq\max\left\{\int_{\Omega}u_{0},\frac{\lambda_{+}|\Omega|}{\mu}\right\}\ \ \ \mbox{for all}\ t>0.
\end{equation}
\item[(2)]
\begin{equation}\label{d1}
\int_{\Omega}u_{\varepsilon}(\cdot,t)\leq\frac{2\hat{\lambda}|\Omega|}{\mu}\ \ \ \mbox{for all}\ t\geq \frac{\ln 2}{\hat\lambda}
\end{equation}
as well as
\begin{equation}\label{d2}
\int_{t}^{t+L}\int_{\Omega}u_{\varepsilon}^{2}
\leq\frac{2|\Omega|}{\mu^{2}}(\hat{\lambda}+\hat{\lambda}^{2}L)\ \ \ \mbox{for all}\ t\geq \frac{\ln 2}{\hat\lambda}\ \mbox{and}\ L>0.
\end{equation}
\item[(3)]
\begin{equation}\label{d3}
\int_{t}^{t+L}\left(\int_{\Omega}v_\ep^2+\int_{\Omega}|\Delta v_\ep|^2\right)\leq\int_{\Omega}v_\ep^2(\cdot,t)+\int_{\Omega}|\nabla v_\ep(\cdot,t)|^2+2\int_{t}^{t+L}\int_{\Omega}u_\ep^2\ \ \ \mbox{for all}\ t>0.
\end{equation}
\item[(4)]
\begin{equation}\label{d4}
\int_{\Omega}v_\ep^2(\cdot, t)+\int_{\Omega}|\nabla v_\ep(\cdot, t)|^2\leq M_{0}\ \ \ \mbox{for all}\ t>0
\end{equation}
with
\begin{equation}\label{d5}
M_{0}:=\max\left\{\int_{\Omega}v_0^2+\int_{\Omega}|\nabla v_{0}|^2+\frac{2}{\mu}\int_{\Omega}u_0,\frac{2(\lambda_{+}+1)}{\mu}\cdot\max\left\{\int_{\Omega}u_{0},\frac{\lambda_{+}|\Omega|}{\mu}\right\}\right\}.
\end{equation}
\end{enumerate}
\end{lemma}

\begin{proof}
Fixing $\bar{\lambda}\geq\lambda_{+}$, we integrate the first equation in \eqref{model111} and employ the H\"older inequality to obtain
\begin{equation}\label{6}
\begin{aligned}
\frac{\rm d}{{\rm d}t}\int_{\Omega}u_{\varepsilon}\leq& \bar{\lambda}\int_{\Omega}u_{\varepsilon}-\mu\int_{\Omega}u_{\varepsilon}^{2}-\varepsilon\int_{\Omega}u_{\varepsilon}^{\theta}\\
\leq&\bar{\lambda}\int_{\Omega}u_{\varepsilon}-\frac{\mu}{|\Omega|}\left(\int_{\Omega}u_{\varepsilon}\right)^{2}\ \ \ \mbox{for all}\ t>0\ \mbox{and}\ \varepsilon\in(0,1),
\end{aligned}
\end{equation}
which leads to \eqref{d0} by an ODE comparison argument and the choice $\bar{\lambda}=\lambda_{+}$. Setting $\bar{\lambda}=\hat{\lambda}>0$ and solving the Bernoulli-type ODI in \eqref{6}, we can estimate
\begin{equation*}
\begin{aligned}
\int_\Omega u_\varepsilon(t)\le& \frac{1}{\frac{e^{-\hat\lambda t}}{\int_\Omega u_{0}}+\frac{\mu}{\hat\lambda|\Omega|}\cdot(1-e^{-\hat\lambda t})}\\
\leq&\frac{1}{\frac{\mu}{\hat\lambda|\Omega|}\cdot(1-e^{-\hat\lambda t})}\ \ \ \mbox{for all}\ t>0\ \mbox{and}\ \varepsilon\in(0,1),
\end{aligned}
\end{equation*}
which immediately gives \eqref{d1}, because $1-e^{-\hat\lambda t}\geq \frac{1}{2}$ for all $t\geq \frac{\ln 2}{\hat\lambda}$. Based on \eqref{d1}, integrating \eqref{6} in time shows that
\begin{equation*}
\begin{aligned}
\int_{\Omega}u_{\varepsilon}(\cdot,t+L)+\mu\int_{t}^{t+L}\int_{\Omega}u_{\varepsilon}^{2}\leq& \int_{\Omega}u_{\varepsilon}(\cdot,t)+\hat\lambda\int_{t}^{t+L}\int_{\Omega}u_{\varepsilon}\\
\leq&\frac{2|\Omega|}{\mu}(\hat{\lambda}+\hat{\lambda}^{2}L)\ \ \ \mbox{for all}\ t\geq \frac{\ln 2}{\hat\lambda}\ \mbox{and}\ \varepsilon\in(0,1),
\end{aligned}
\end{equation*}
and hence \eqref{d2} follows. Then multiplying the second equation in \eqref{model111} by $2v_\ep$ and $-2\Delta v_\ep$, respectively, we deduce that
\begin{equation}\label{1666}
\frac{{\rm d}}{{\rm d}t}\left(\int_{\Omega}v_\ep^2+\int_{\Omega}|\nabla v_\ep|^2\right)
+4\int_{\Omega}|\nabla v_\ep|^2+2\int_{\Omega}|\Delta v_\ep|^2+2\int_{\Omega}v_\ep^2=2\int_{\Omega}u_\ep v_\ep-2\int_{\Omega}u_\ep\Delta v_\ep,
\end{equation}
where we may use Young's inequality to derive
\begin{equation*}
2\int_{\Omega}u_\ep v_\ep-2\int_{\Omega}u_\ep\Delta v_\ep\leq \int_{\Omega}v_\ep^{2}+\int_{\Omega}|\Delta v_\ep|^2+2\int_{\Omega}u_\ep^{2}\ \ \ \mbox{for all}\ t>0\ \mbox{and}\ \varepsilon\in(0,1).
\end{equation*}
Therefore, \eqref{1666} entails that
\begin{equation}\label{1777}
\frac{{\rm d}}{{\rm d}t}\left(\int_{\Omega}v_\ep^2+\int_{\Omega}|\nabla v_\ep|^2\right)
+4\int_{\Omega}|\nabla v_\ep|^2+\int_{\Omega}|\Delta v_\ep|^2+\int_{\Omega}v_\ep^2\leq 2\int_{\Omega}u_\ep^{2}
\end{equation}
for all $t>0$ and $\varepsilon\in(0,1)$, so that \eqref{d3} can be concluded upon integrating in time. Moreover, we combine \eqref{d0} and \eqref{6} with \eqref{1777} to see that
\begin{equation*}
y_{\varepsilon}(t):=\int_{\Omega}v_\ep^2(\cdot, t)+\int_{\Omega}|\nabla v_\ep(\cdot, t)|^2+\frac{2}{\mu}\int_{\Omega}u_\ep(\cdot, t)
\end{equation*}
fulfills
\begin{equation*}
\begin{aligned}
y'_{\varepsilon}(t)+y_{\varepsilon}(t)\leq& \frac{2(\lambda_{+}+1)}{\mu}\int_{\Omega}u_{\varepsilon}\\
\leq&\frac{2(\lambda_{+}+1)}{\mu}\cdot\max\left\{\int_{\Omega}u_{0},\frac{\lambda_{+}|\Omega|}{\mu}\right\}\ \ \ \mbox{for all}\ t>0\ \mbox{and}\ \varepsilon\in(0,1),
\end{aligned}
\end{equation*}
which readily yields \eqref{d4}-\eqref{d5}.
\end{proof}

When $N=4$, the following interpolation inequality plays a key role in controlling $\int_\Omega |\nabla v_\varepsilon|^6$.
\begin{lemma}\label{a}
Let $N=4$. Then there exists $C_{\star}(\Omega)>0$ such that for all $t>0$ and $\varepsilon\in(0,1)$, we have
\begin{equation}\label{a1}
\begin{aligned}
\int_{\Omega}|\nabla v_{\varepsilon}|^6
\le C_{\star}\left(\int_{\Omega}|\nabla v_{\varepsilon}|^4\right)^{\frac{1}{2}}\left(\int_{\Omega}\left|\nabla|\nabla v_{\varepsilon}|^2\right|^2+\int_{\Omega}|\nabla v_{\varepsilon}|^4\right).
\end{aligned}
\end{equation}
\end{lemma}
\begin{proof}
Let $w_{\varepsilon}=|\nabla v_{\varepsilon}|^2$. It follows from the H\"older inequality that
\begin{equation}\label{1}
\int_{\Omega}w_{\varepsilon}^{3}\leq \left(\int_{\Omega}w_{\varepsilon}^{2}\right)^{\frac{1}{2}}\left(\int_{\Omega}w_{\varepsilon}^{4}\right)^{\frac{1}{2}}\ \ \ \mbox{for all}\ t>0\ \mbox{and}\ \varepsilon\in(0,1).
\end{equation}
Here the Sobolev embedding $W^{1,2}(\Omega)\hookrightarrow L^{4}(\Omega)$ provides $C_{\star}(\Omega)>0$ such that
\begin{equation}\label{2}
\|w_{\varepsilon}\|_{L^{4}(\Omega)}^{2}\leq C_{\star}(\Omega)\left(\|w_{\varepsilon}\|_{L^{2}(\Omega)}^{2}+\|\nabla w_{\varepsilon}\|_{L^{2}(\Omega)}^{2}\right).
\end{equation}
Thus the combination of \eqref{1} and \eqref{2} gives \eqref{a1}.
\end{proof}

We now combine the preceding estimates with the key interpolation inequality \eqref{a1} to derive the following energy inequality.
\begin{lemma}\label{b}
Let $N=4$. For each $\varepsilon\in(0,1)$, one has
\begin{equation}\label{b1}
\begin{aligned}
&\frac{{\rm d}}{{\rm d}t}\left(\int_{\Omega}u_\ep^2+\int_{\Omega}|\nabla v_\ep|^4\right)+\int_{\Omega}|\nabla u_\ep|^2
+2\int_{\Omega}\left|\nabla|\nabla v_\ep|^2\right|^2+4\int_{\Omega}|\nabla v_\ep|^4\\
\le& 2\lambda\int_{\Omega}u_\ep^2+C_{\star}\left(8+\frac{8\chi^{6}}{27\mu^2}\right)\left(\int_{\Omega}|\nabla v_{\varepsilon}|^4\right)^{\frac{1}{2}}\left(\int_{\Omega}\left|\nabla|\nabla v_{\varepsilon}|^2\right|^2+\int_{\Omega}|\nabla v_{\varepsilon}|^4\right)\ \ \ \mbox{for all}\ t>0.
\end{aligned}
\end{equation}
\end{lemma}
\begin{proof}
Testing the first equation in \eqref{model111} by $2u_\ep$ and integrating by parts, we see that
\begin{equation}\label{3}
\frac{{\rm d}}{{\rm d}t}\int_{\Omega}u_\ep^2+2\int_{\Omega}|\nabla u_\ep|^2 \leq2\chi\int_{\Omega}u_\ep\nabla u_\ep\cdot\nabla v_\ep+2\lambda\int_{\Omega}u_\ep^2-2\mu\int_{\Omega}u_\ep^3-2\ep\int_{\Omega}u_\ep^{\theta+1}
\end{equation}
for all $t>0$ and $\varepsilon\in(0,1)$, where two applications of Young's inequality show that
\begin{equation}\label{4}
\begin{aligned}
2\chi\int_{\Omega}u_\ep\nabla u_\ep\cdot\nabla v_\ep\leq& \frac12\int_{\Omega}|\nabla u_\ep|^2+2\chi^{2}\int_{\Omega}u_\ep^2|\nabla v_\ep|^2\\
\leq&\frac12\int_{\Omega}|\nabla u_\ep|^2+2\mu\int_{\Omega}u_\ep^3+\frac{8\chi^{6}}{27\mu^2}\int_{\Omega}|\nabla v_\ep|^6.
\end{aligned}
\end{equation}
Combining \eqref{4} with \eqref{3} implies
\begin{equation}\label{5}
\frac{{\rm d}}{{\rm d}t}\int_{\Omega}u_\ep^2+\frac{3}{2}\int_{\Omega}|\nabla u_\ep|^2 \leq 2\lambda\int_{\Omega}u_\ep^2+\frac{8\chi^{6}}{27\mu^2}\int_{\Omega}|\nabla v_\ep|^6\ \ \ \mbox{for all}\ t>0\ \mbox{and}\ \varepsilon\in(0,1).
\end{equation}
According to \cite[Lemma 4.6]{L2015}, we further use Young's inequality to estimate that
\begin{equation*}
\begin{aligned}
\frac{{\rm d}}{{\rm d}t}\int_{\Omega}|\nabla v_\ep|^4+2\int_{\Omega}\left|\nabla|\nabla v_\ep|^2\right|^2+4\int_{\Omega}|\nabla v_\ep|^4\leq& 4\int_{\Omega}|\nabla v_\ep|^3|\nabla u_\ep|\\
\leq&\frac12\int_{\Omega}|\nabla u_\ep|^2+8\int_{\Omega}|\nabla v_\ep|^6
\end{aligned}
\end{equation*}
for all $t>0$ and $\varepsilon\in(0,1)$, which together with \eqref{5} and \eqref{a1} proves \eqref{b1}.
\end{proof}

The next lemma shows that for suitably small growth, once the coupled quality becomes sufficiently small
at some large time, it remains small thereafter.
\begin{lemma}\label{c}
Let $N=4$ and
\begin{equation}\label{c1}
0<\eta< \min\left\{1, \frac{1}{C^{2}_\star\left(8+\frac{8\chi^{6}}{27\mu^2}\right)^{2}}\right\}.
\end{equation}
Suppose that $\lambda_+<\hat\lambda$ with
\begin{equation}\label{c2}
0<\hat{\lambda}\le \min\left\{\frac1{4C_P}, \sqrt{\frac{\eta C_P\mu^2\cdot\min\left\{\frac1{2C_P},3\right\}}{16|\Omega|}}\right\}.
\end{equation}
For every $\varepsilon\in(0,1)$, if
\begin{equation}\label{c4}
\int_{\Omega}u_\ep^2(\cdot,T_{\star})+\int_{\Omega}|\nabla v_\ep(\cdot,T_{\star})|^4<\frac{\eta}{2} \ \ \ \mbox{for some}\ T_{\star}> \frac{\ln 2}{\hat\lambda},
\end{equation}
then one has
\begin{equation}\label{c3}
\int_{\Omega}u_\ep^2(\cdot,t) +\int_{\Omega}|\nabla v_\ep(\cdot,t)|^4<\eta\ \ \ \mbox{for all}\ t\geq T_{\star}.
\end{equation}

\end{lemma}
\begin{proof}
In view of \eqref{d1} and Lemma \ref{p1}, we arrive at
\begin{equation*}
\begin{aligned}
-\int_{\Omega}|\nabla u_{\varepsilon}|^{2}\leq&-\frac{1}{C_{P}}\int_{\Omega}u_{\varepsilon}^{2}+\frac{1}{C_{P}|\Omega|}\left(\int_{\Omega}u_{\varepsilon}\right)^{2}\\
\leq&-\frac{1}{C_{P}}\int_{\Omega}u_{\varepsilon}^{2}+\frac{4|\Omega|\hat{\lambda}^{2}}{C_{P}\mu^{2}}\ \ \ \mbox{for all}\ t\geq \frac{\ln 2}{\hat\lambda}\ \mbox{and}\ \varepsilon\in(0,1),
\end{aligned}
\end{equation*}
which updates \eqref{b1} as
\begin{equation}\label{11}
\begin{aligned}
&\frac{{\rm d}}{{\rm d}t}\left(\int_{\Omega}u_\ep^2+\int_{\Omega}|\nabla v_\ep|^4\right)
+2\int_{\Omega}\left|\nabla|\nabla v_\ep|^2\right|^2+4\int_{\Omega}|\nabla v_\ep|^4+\left(\frac1{C_P}-2\hat{\lambda}\right)\int_{\Omega}u_\ep^2\\
\le& C_{\star}\left(8+\frac{8\chi^{6}}{27\mu^2}\right)\left(\int_{\Omega}|\nabla v_{\varepsilon}|^4\right)^{\frac{1}{2}}\left(\int_{\Omega}\left|\nabla|\nabla v_{\varepsilon}|^2\right|^2+\int_{\Omega}|\nabla v_{\varepsilon}|^4\right)+\frac{4|\Omega|\hat{\lambda}^{2}}{C_{P}\mu^{2}}
\end{aligned}
\end{equation}
for all $t\geq \frac{\ln 2}{\hat\lambda}$ and $\varepsilon\in(0,1)$. Letting $J_{\varepsilon}(t):=\int_{\Omega}u_\ep^2(\cdot, t)+\int_{\Omega}|\nabla v_\ep(\cdot, t)|^4$ with $t\geq T_{\star}>\frac{\ln 2}{\hat\lambda}$, we find that for each $\varepsilon\in(0,1)$,
\begin{equation}\label{12}
T_{\varepsilon}:=\sup\left\{T'>T_{\star}\;\bigg|\;J_{\varepsilon}(t)\leq\eta\ \ \ \mbox{for all} \ t\in[T_{\star},T')\right\}
\end{equation}
is a well-defined  element of $(0,\infty]$ thanks to the continuity of $J_{\varepsilon}(t)$ and \eqref{c4}.
Then it follows from \eqref{11} and \eqref{12} that
\begin{equation*}
\begin{aligned}
&\frac{{\rm d}}{{\rm d}t}\left(\int_{\Omega}u_\ep^2+\int_{\Omega}|\nabla v_\ep|^4\right)
+2\int_{\Omega}\left|\nabla|\nabla v_\ep|^2\right|^2+4\int_{\Omega}|\nabla v_\ep|^4+\left(\frac1{C_P}-2\hat{\lambda}\right)\int_{\Omega}u_\ep^2\\
\le& C_{\star}\left(8+\frac{8\chi^{6}}{27\mu^2}\right)\eta^{\frac{1}{2}}\left(\int_{\Omega}\left|\nabla|\nabla v_{\varepsilon}|^2\right|^2+\int_{\Omega}|\nabla v_{\varepsilon}|^4\right)+\frac{4|\Omega|\hat{\lambda}^{2}}{C_{P}\mu^{2}}\ \ \ \mbox{for all}\ t\in[T_{\star},T_{\varepsilon}),
\end{aligned}
\end{equation*}
which implies
\begin{equation}\label{13}
J'_{\varepsilon}(t)+\min\left\{\frac1{2C_P},3\right\}
\cdot J_{\varepsilon}(t)\le\frac{4|\Omega|\hat{\lambda}^{2}}{C_{P}\mu^{2}}\ \ \ \mbox{for all}\ t\in[T_{\star},T_{\varepsilon})
\end{equation}
because $-\left(\frac1{C_P}-2\hat{\lambda}\right)\leq -\frac1{2C_P}$ and $C_{\star}\left(8+\frac{8\chi^{6}}{27\mu^2}\right)\eta^{\frac{1}{2}}\leq 1$ by \eqref{c1}-\eqref{c2}. Solving the ODI \eqref{13} then yields
\begin{equation}\label{14}
\begin{aligned}
J_{\varepsilon}(t)
\le&
e^{-\min\left\{\frac1{2C_P},3\right\}(t-T_{\star})}J_{\varepsilon}(T_{\star})
+\frac{4|\Omega|\hat{\lambda}^{2}}{\min\left\{\frac1{2C_P},3\right\}\cdot C_{P}\mu^{2}}
\left(1-e^{-\min\left\{\frac1{2C_P},3\right\}(t-T_{\star})}\right)\\
\le&e^{-\min\left\{\frac1{2C_P},3\right\}(t-T_{\star})}J_{\varepsilon}(T_{\star})+\frac{\eta}{4}\left(1-e^{-\min\left\{\frac1{2C_P},3\right\}(t-T_{\star})}\right)\\
\le&\max\left\{J_{\varepsilon}(T_{\star}), \frac{\eta}{4}\right\}\leq\frac{\eta}{2}\ \ \ \mbox{for all}\ t\in[T_{\star},T_{\varepsilon}).
\end{aligned}
\end{equation}
By the continuity of $J_{\varepsilon}(t)$, the estimate \eqref{14} rules out the possibility that $T_{\varepsilon}<\infty$, and hence establishes \eqref{c3}.
\end{proof}

It remains to verify that the smallness condition in
Lemma~\ref{c} is attained, thereby establishing the
eventual smallness of $\|u_{\varepsilon}(\cdot,t)\|_{L^{2}(\Omega)}$ and $\|\nabla v_{\varepsilon}(\cdot,t)\|_{L^{4}(\Omega)}$.
\begin{proposition}\label{e}
Let $N=4$. For every $\delta>0$, there exists $\lambda_{\delta,1}=\lambda_{\delta,1}(\delta,\chi,\mu,\Omega)$ such that whenever $\lambda<\lambda_{\delta,1}$, one can find $T_{\delta,1}=T_{\delta,1}(\delta,\chi,\mu,\Omega,u_{0},v_{0})>0$ such that
\begin{equation}\label{e1}
\int_{\Omega}u_\ep^2(\cdot,t)+\int_{\Omega}|\nabla v_\ep(\cdot,t)|^4<\delta
\ \ \ \mbox{for all}\ t\geq T_{\delta,1}\ \mbox{and}\ \varepsilon\in(0,1).
\end{equation}
\end{proposition}

\begin{proof}
For any $\delta>0$, we define
\begin{equation}\label{eta}
\eta_{\delta}:= \frac{1}{2}\min\left\{\delta, 1, \frac{1}{C^{2}_\star\left(8+\frac{8\chi^{6}}{27\mu^2}\right)^{2}}\right\}.
\end{equation}
Let $C_\Omega>0$ be the constant in \eqref{17} below. Choose $\lambda_{\delta,1}>0$ sufficiently small so that \eqref{c2} holds with $\hat{\lambda}=\lambda_{\delta,1}$ and, moreover,
\begin{equation}\label{15}
\frac{2(1+2C_\Omega)|\Omega|}{\mu^2}\lambda_{\delta,1}^2<\frac{\eta_{\delta}}{8}.
\end{equation}
Since the right-hand side of \eqref{d5} is nondecreasing in
$\lambda_+$, the bound $\lambda_+<\lambda_{\delta,1}$ yields a
uniform upper bound for $M_0$. We may therefore choose $L>0$
suitably large that
\begin{equation}\label{16}
\frac{1}{L}\left(\frac{2(1+2C_\Omega)|\Omega|}{\mu^2}\lambda_{\delta,1}+C_\Omega M_{0}\right)<\frac{\eta_{\delta}}{8}
\end{equation}
holds for every $\lambda<\lambda_{\delta,1}$.
As the Sobolev embedding $W^{2,2}(\Omega)\hookrightarrow W^{1,4}(\Omega)$, together with the Neumann elliptic regularity warrants
\begin{equation}\label{17}
\left(\int_{\Omega}|\nabla v_{\varepsilon}|^{4}\right)^{\frac{1}{2}}\leq C_{\Omega}\left(\int_{\Omega}v_{\varepsilon}^{2}+\int_{\Omega}|\Delta v_{\varepsilon}|^{2}\right),
\end{equation}
we use \eqref{d2}-\eqref{d5} to infer that for $T_{\star}=\frac{\ln 2}{\lambda_{\delta,1}}+1$ and any $L>0$,
\begin{equation*}
\begin{aligned}
&\int_{T_{\star}}^{T_{\star}+L}\left\{\int_{\Omega}u_\ep^2
+\left(\int_{\Omega}|\nabla v_\ep|^4\right)^{1/2}\right\}\\
\le&\int_{T_{\star}}^{T_{\star}+L}\int_{\Omega}u_\ep^2+C_\Omega\int_{T_{\star}}^{T_{\star}+L}\left(\int_{\Omega}v_\ep^2+\int_{\Omega}|\Delta v_\ep|^2\right)\\
\le&(1+2C_\Omega)\int_{T_{\star}}^{T_{\star}+L}\int_{\Omega}u_\ep^2+C_\Omega\left(\int_{\Omega}v_\ep^2(\cdot,T_{\star})+\int_{\Omega}|\nabla v_\ep(\cdot,T_{\star})|^2\right)\\
\le&\frac{2(1+2C_\Omega)|\Omega|}{\mu^2}\lambda_{\delta,1}^2L+\frac{2(1+2C_\Omega)|\Omega|}{\mu^2}\lambda_{\delta,1}+C_\Omega M_{0}\ \ \ \mbox{for all}\ \varepsilon\in(0,1),
\end{aligned}
\end{equation*}
which further implies
\begin{equation}\label{18}
\begin{aligned}
\frac{1}{L}\int_{T_{\star}}^{T_{\star}+L}\left\{\int_{\Omega}u_\ep^2
+\left(\int_{\Omega}|\nabla v_\ep|^4\right)^{\frac{1}{2}}\right\}\le&\frac{2(1+2C_\Omega)|\Omega|}{\mu^2}\lambda_{\delta,1}^2+\frac{1}{L}\left(\frac{2(1+2C_\Omega)|\Omega|}{\mu^2}\lambda_{\delta,1}+C_\Omega M_{0}\right)\\
\le&\frac{\eta_{\delta}}{8}+\frac{\eta_{\delta}}{8}=\frac{\eta_{\delta}}{4}\ \ \ \mbox{for all}\ \varepsilon\in(0,1)
\end{aligned}
\end{equation}
thanks to \eqref{15} and \eqref{16}. On the basis of \eqref{18}, for each $\varepsilon\in(0,1)$ we can pick $\widetilde{T_{\varepsilon}} \in (T_{\star},T_{\star}+L)$ such that
\begin{equation*}
\int_{\Omega}u_\ep^2(\cdot, \widetilde{T_{\varepsilon}})+\left(\int_{\Omega}|\nabla v_\ep(\cdot, \widetilde{T_{\varepsilon}})|^4\right)^{1/2}\leq\frac{\eta_{\delta}}{4},
\end{equation*}
and hence,
\begin{equation*}
\int_{\Omega}u_\ep^2(\cdot, \widetilde{T_{\varepsilon}})+\int_{\Omega}|\nabla v_\ep(\cdot, \widetilde{T_{\varepsilon}})|^4\leq\frac{\eta_{\delta}}{4}+\frac{\eta_{\delta}^{2}}{16}<\frac{\eta_{\delta}}{2}.
\end{equation*}
Then an application of \eqref{eta} and Lemma \ref{c} ensures
\begin{equation}\label{20}
\int_{\Omega}u_\ep^2(\cdot, t)+\int_{\Omega}|\nabla v_\ep(\cdot, t)|^4\leq\eta_{\delta}<\delta\ \ \ \mbox{for all}\ t\geq \widetilde{T_{\varepsilon}}.
\end{equation}
Since $\widetilde{T_{\varepsilon}}<T_{\star}+L$, \eqref{e1} becomes a consequence of \eqref{20} upon setting $T_{\delta,1}:=T_{\star}+L$.
\end{proof}

With this eventual smallness available, we use semigroup estimates to improve the eventual smallness of $u_\varepsilon$ in $L^3(\Omega)$ and $\nabla v_\varepsilon$ in $L^6(\Omega)$.
\begin{lemma}\label{f}
Let $N=4$. There exist $C>0$ and $\delta_\star>0$ such that for any $\delta\in(0, \delta_\star]$, one can find $\lambda_{\delta,2}=\lambda_{\delta,2}(\delta,\chi,\mu,\Omega)$ with the property that if $\lambda<\lambda_{\delta,2}$, then there exists $T_{\delta,2}=T_{\delta,2}(\delta,\chi,\mu,\Omega,u_{0},v_{0})>T_{\delta,1}$ so that for any choice of $\varepsilon\in(0,1)$, we have
\begin{equation}\label{f1}
\|u_{\varepsilon}(\cdot,t)\|_{L^{3}(\Omega)}\leq C\delta^{\frac{1}{2}}\ \ \ \mbox{for all} \  t\geq T_{\delta,2}
\end{equation}
and
\begin{equation}\label{f2}
\|\nabla v_{\varepsilon}(\cdot,t)\|_{L^{6}(\Omega)}\leq C\delta^{\frac{1}{4}}\ \ \ \mbox{for all} \  t\geq T_{\delta,2}.
\end{equation}
\end{lemma}
\begin{proof}
The regularization estimates for the Neumann heat semigroup $\{e^{t\Delta}\}_{t\geq0}$ provide some positive constants $C_{i}$ $(i=1, ..., 5)$ such that
\begin{equation}\label{23}
\| e^{\tau\Delta}\psi \|_{L^3(\Omega)}\leq C_{1}\cdot\left(1+\tau^{-\frac{1}{3}}\right)\|\psi\|_{L^2(\Omega)}\ \ \ \mbox{for all}\ \psi\in L^2(\Omega)\ \mbox{and}\ \tau>0
\end{equation}
and
\begin{equation}\label{24}
\| e^{\tau\Delta}\psi \|_{L^3(\Omega)}\leq C_{2}\cdot\|\psi\|_{L^3(\Omega)}\ \ \ \mbox{for all}\ \psi\in L^3(\Omega)\ \mbox{and}\ \tau>0,
\end{equation}
and that
\begin{equation}\label{25}
\| e^{\tau\Delta}\nabla\cdot\psi \|_{L^3(\Omega)}\leq C_{3}\cdot\left(1+\tau^{-\frac{5}{6}}\right)\|\psi\|_{L^2(\Omega)}\ \ \ \mbox{for all}\ \psi\in (L^{2}(\Omega))^{N}\ \mbox{and}\ \tau>0,
\end{equation}
and, again,
\begin{equation}\label{26}
\|\nabla e^{\tau\Delta}\psi \|_{L^6(\Omega)}\leq C_{4}\cdot\left(1+\tau^{-\frac{1}{6}}\right)\|\nabla\psi\|_{L^4(\Omega)}\ \ \ \mbox{for all}\ \psi\in W^{1,6}(\Omega)\ \mbox{and}\ \tau>0
\end{equation}
as well as
\begin{equation}\label{27}
\|\nabla e^{\tau\Delta}\psi \|_{L^6(\Omega)}\leq C_{5}\cdot\left(1+\tau^{-\frac{5}{6}}\right)\|\psi\|_{L^3(\Omega)}\ \ \ \mbox{for all}\ \psi\in L^{3}(\Omega)\ \mbox{and}\ \tau>0.
\end{equation}
Moreover, we may choose $C_{6}>0$ and $C_{7}>0$ such that
\begin{equation}\label{28}
\sigma^{\frac{1}{3}}\int_{0}^{\sigma}s^{-\frac{1}{2}}\left(1+(\sigma-s)^{-\frac{5}{6}}\right)ds\leq C_{6}\ \ \ \mbox{for all}\ \sigma\in(0,2]
\end{equation}
and
\begin{equation}\label{29}
\sigma^{\frac{1}{6}}\int_{0}^{\sigma} s^{-\frac{1}{3}}\left(1+(\sigma-s)^{-\frac{5}{6}}\right)ds\leq C_{7}\ \ \ \mbox{for all}\ \sigma\in(0,2].
\end{equation}
We then abbreviate
\begin{equation}\label{21}
\delta_{\star}=\min\left\{\frac{1}{(48\chi C_{3}C_{4}C_{6})^{4}}, \left(\frac{C_{4}}{2C_{1}C_{5}C_{7}}\right)^{4}\right\}
\end{equation}
and fix $\delta\in(0,\delta_\star]$. Let $\lambda_{\delta,1}>0$ denote the threshold provided by Proposition \ref{e}, and set
\begin{equation}\label{22}
\lambda_{\delta,2}:=\min\left\{\lambda_{\delta,1},\frac{1}{12C_{2}}\right\}>0.
\end{equation}
Assume that $\lambda\leq\lambda_{\delta,2}$. Since $\lambda_{\delta,2}<\lambda_{\delta,1}$ by \eqref{22}, Proposition \ref{e} ensures the existence of $T_{\delta,1}(\delta,\chi,\mu,\Omega,u_{0},v_{0})>0$ such that, upon fixing any $t_0\geq T_{\delta,2}-2:=T_{\delta,1}$, one has
\begin{equation}\label{30}
\int_\Omega u_\varepsilon^2(\cdot,t_{0})+\int_\Omega|\nabla v_\varepsilon(\cdot,t_{0})|^4 <\delta \ \ \ \mbox{for all} \  \varepsilon\in(0,1).
\end{equation}
For each $\varepsilon\in(0,1)$ and $\sigma\in(0,2]$, we define
\begin{equation}\label{31}
X_{\varepsilon}(\sigma):=\|\sigma^{\frac{1}{3}}u_{\varepsilon}(\cdot,t_{0}+\sigma)\|_{L^{3}(\Omega)}
\end{equation}
as well as
\begin{equation}\label{32}
Y_{\varepsilon}(\sigma):=\|\sigma^{\frac{1}{6}}\nabla v_{\varepsilon}(\cdot,t_{0}+\sigma)\|_{L^{6}(\Omega)}.
\end{equation}
By the continuity of $u_{\varepsilon}$ and $\nabla v_{\varepsilon}$, we infer that
\begin{equation}\label{33}
S_{\varepsilon}:=\left\{T'\in(0,2]\;\bigg|\;X_{\varepsilon}(\sigma)< 12 C_{1}\delta^{\frac{1}{2}}\ \  \mbox{and} \ \ Y_{\varepsilon}(\sigma)< 12 C_{4}\delta^{\frac{1}{4}}\ \ \ \mbox{for all} \ \sigma\in(0, T')\right\}
\end{equation}
is not empty, and hence $T_\varepsilon:=\sup S_{\varepsilon}$ is a well-defined element of $(0,2]$. We next apply the variation-of-constants formula and the order-preserving of $\{e^{\tau\Delta}\}_{\tau\geq0}$ to derive
\begin{equation}\label{341}
\begin{aligned}
&u_{\varepsilon}(\cdot,t_{0}+\sigma)=e^{\sigma\Delta}u_{\varepsilon}(\cdot,t_{0})-\chi\int_{0}^{\sigma}e^{(\sigma-s)\Delta}\nabla \cdot\left(u_{\varepsilon}(\cdot,t_{0}+s)\nabla v_{\varepsilon}(\cdot,t_{0}+s)\right)ds\\
&+\int_{0}^{\sigma}e^{(\sigma-s)\Delta}\left(\lambda u_{\varepsilon}(\cdot,t_{0}+s)-\mu u_{\varepsilon}^{2}(\cdot,t_{0}+s)-\varepsilon u^{\theta}_{\varepsilon}(\cdot,t_{0}+s)\right)ds\\
\leq&e^{\sigma\Delta}u_{\varepsilon}(\cdot,t_{0})-\chi\int_{0}^{\sigma}e^{(\sigma-s)\Delta}\nabla \cdot\left(u_{\varepsilon}(\cdot,t_{0}+s)\nabla v_{\varepsilon}(\cdot,t_{0}+s)\right)ds+\lambda_{+}\int_{0}^{\sigma}e^{(\sigma-s)\Delta}u_{\varepsilon}(\cdot,t_{0}+s)ds\\
\end{aligned}
\end{equation}
and
\begin{equation*}
\nabla v_{\varepsilon}(\cdot,t_{0}+\sigma)=\nabla e^{\sigma(\Delta-1)}v_{\varepsilon}(\cdot,t_{0})+\int_{0}^{\sigma}\nabla e^{(\sigma-s)(\Delta-1)}u_{\varepsilon}(\cdot,t_{0}+s)ds
\end{equation*}
for all $\sigma\in(0,T_{\varepsilon})$, which yields
\begin{equation}\label{34}
\begin{aligned}
\|\sigma^{\frac{1}{3}}u_{\varepsilon}(\cdot,t_{0}+\sigma)\|_{L^{3}(\Omega)}\leq&\underbrace{\sigma^{\frac{1}{3}}\|e^{\sigma\Delta}u_{\varepsilon}(\cdot,t_{0})\|_{L^{3}(\Omega)}}\limits_{:=I_{1}}+\underbrace{\lambda_{+}\sigma^{\frac{1}{3}}\int_{0}^{\sigma}\|e^{(\sigma-s)\Delta}u_{\varepsilon}(\cdot,t_{0}+s)\|_{L^{3}(\Omega)}ds}\limits_{:=I_{2}}\\
&+\underbrace{\chi\sigma^{\frac{1}{3}}\int_{0}^{\sigma}\|e^{(\sigma-s)\Delta}\nabla \cdot\left(u_{\varepsilon}(\cdot,t_{0}+s)\nabla v_{\varepsilon}(\cdot,t_{0}+s)\right)\|_{L^{3}(\Omega)}ds}\limits_{:=I_{3}}\\
\end{aligned}
\end{equation}
and
\begin{equation}\label{35}
\begin{aligned}
&\|\sigma^{\frac{1}{6}}\nabla v_{\varepsilon}(\cdot,t_{0}+\sigma)\|_{L^{6}(\Omega)}\\
\leq&\underbrace{\sigma^{\frac{1}{6}}\|\nabla e^{\sigma(\Delta-1)} v_{\varepsilon}(\cdot,t_{0})\|_{L^{6}(\Omega)}}\limits_{:=I_{4}}+\underbrace{\sigma^{\frac{1}{6}}\int_{0}^{\sigma}\|\nabla e^{(\sigma-s)(\Delta-1)}u_{\varepsilon}(\cdot,t_{0}+s)\|_{L^{6}(\Omega)}ds}\limits_{:=I_{5}}
\end{aligned}
\end{equation}
for all $\sigma\in(0,T_{\varepsilon})$. Here, it follows from \eqref{23}, \eqref{26} and \eqref{30} that
\begin{equation}\label{36}
\begin{aligned}
I_{1}\leq& C_{1}\sigma^{\frac{1}{3}}\left(1+\sigma^{-\frac{1}{3}}\right)\|u_{\varepsilon}(\cdot,t_{0})\|_{L^{2}(\Omega)}\\
\leq& C_{1}\left(1+\sigma^{\frac{1}{3}}\right)\delta^{\frac{1}{2}}\leq 3C_{1}\delta^{\frac{1}{2}}\ \ \ \mbox{for all} \  \sigma\in(0,T_{\varepsilon})
\end{aligned}
\end{equation}
and
\begin{equation}\label{37}
\begin{aligned}
I_{4}\leq& C_{4}\sigma^{\frac{1}{6}}\left(1+\sigma^{-\frac{1}{6}}\right)\|\nabla v_{\varepsilon}(\cdot,t_{0})\|_{L^{4}(\Omega)}\\
\leq& C_{4}\left(1+\sigma^{\frac{1}{6}}\right)\delta^{\frac{1}{4}}\leq 3C_{4}\delta^{\frac{1}{4}}\ \ \ \mbox{for all} \  \sigma\in(0,T_{\varepsilon}).
\end{aligned}
\end{equation}
Next, based on the definitions in \eqref{31}-\eqref{33}, we may employ \eqref{24} and \eqref{27} together with \eqref{29} to estimate
\begin{equation}\label{38}
\begin{aligned}
I_{2}\leq& C_{2}\lambda_{+}\sigma^{\frac{1}{3}}\int_{0}^{\sigma}\|u_{\varepsilon}(\cdot,t_{0}+s)\|_{L^{3}(\Omega)}ds\\
\leq& 12 C_{1}C_{2}\lambda_{+}\sigma^{\frac{1}{3}}\delta^{\frac{1}{2}}\int_{0}^{\sigma}s^{-\frac{1}{3}}ds\leq 36C_{1}C_{2}\lambda_{+}\delta^{\frac{1}{2}}\ \ \ \mbox{for all} \  \sigma\in(0,T_{\varepsilon})
\end{aligned}
\end{equation}
as well as
\begin{equation}\label{39}
\begin{aligned}
I_{5}\leq& C_{5}\sigma^{\frac{1}{6}}\int_{0}^{\sigma}\left(1+(\sigma-s)^{-\frac{5}{6}}\right)\|u_{\varepsilon}(\cdot,t_{0}+s) \|_{L^{3}(\Omega)}ds\\
\leq& 12 C_{1}C_{5}\delta^{\frac{1}{2}}\cdot\left\{\sigma^{\frac{1}{6}}\int_{0}^{\sigma} s^{-\frac{1}{3}}\left(1+(\sigma-s)^{-\frac{5}{6}}\right)ds\right\}\\
\leq& 12 C_{1}C_{5}C_{7}\delta^{\frac{1}{2}}\ \ \ \mbox{for all} \  \sigma\in(0,T_{\varepsilon}).
\end{aligned}
\end{equation}
Finally, again drawing on \eqref{31}-\eqref{33}, by means of the H\"{o}lder inequality, \eqref{25} and \eqref{28} we obtain
\begin{equation}\label{40}
\begin{aligned}
I_{3}\leq& C_{3}\chi\sigma^{\frac{1}{3}}\int_{0}^{\sigma}\left(1+(\sigma-s)^{-\frac{5}{6}}\right)\|u_{\varepsilon}(\cdot,t_{0}+s) \nabla v_{\varepsilon}(\cdot,t_{0}+s)\|_{L^{2}(\Omega)}ds\\
\leq& C_{3}\chi\sigma^{\frac{1}{3}}\int_{0}^{\sigma}\left(1+(\sigma-s)^{-\frac{5}{6}}\right)\|u_{\varepsilon}(\cdot,t_{0}+s)\|_{L^{3}(\Omega)}\|\nabla v_{\varepsilon}(\cdot,t_{0}+s)\|_{L^{6}(\Omega)}ds\\
\leq& 144\chi C_{1}C_{3}C_{4}\delta^{\frac{3}{4}}\cdot\left\{\sigma^{\frac{1}{3}}\int_{0}^{\sigma}s^{-\frac{1}{2}}\left(1+(\sigma-s)^{-\frac{5}{6}}\right)ds\right\}\\
\leq& 144\chi C_{1}C_{3}C_{4}C_{6}\delta^{\frac{3}{4}}\ \ \ \mbox{for all} \  \sigma\in(0,T_{\varepsilon}).
\end{aligned}
\end{equation}
Now combining \eqref{36}-\eqref{40} with \eqref{34} and \eqref{35} implies
\begin{equation}\label{41}
\begin{aligned}
X_{\varepsilon}(\sigma)\leq& \left(3C_{1}+36C_{1}C_{2}\lambda_{+}+144\chi C_{1}C_{3}C_{4}C_{6}\delta^{\frac{1}{4}}\right)\delta^{\frac{1}{2}}\\
\leq& \left(3C_{1}+3C_{1}+3C_{1}\right)\delta^{\frac{1}{2}}=9C_{1}\delta^{\frac{1}{2}}<12 C_{1}\delta^{\frac{1}{2}}\ \ \ \mbox{for all} \  \sigma\in(0,T_{\varepsilon})
\end{aligned}
\end{equation}
and
\begin{equation}\label{42}
\begin{aligned}
Y_{\varepsilon}(\sigma)\leq& \left(3C_{4}+12 C_{1}C_{5}C_{7}\delta^{\frac{1}{4}}\right)\delta^{\frac{1}{4}}\\
\leq& \left(3C_{4}+6C_{4}\right)\delta^{\frac{1}{4}}=9C_{4}\delta^{\frac{1}{4}}<12 C_{4}\delta^{\frac{1}{4}}\ \ \ \mbox{for all} \  \sigma\in(0,T_{\varepsilon})
\end{aligned}
\end{equation}
due to \eqref{21}-\eqref{22} and the fact that $\lambda_{+}<\lambda_{\delta,2}$. Therefore, in light of the continuity of $u_\varepsilon$ and $\nabla v_{\varepsilon}$, we must have $T_\varepsilon=2$ for all $\varepsilon\in(0,1)$. Consequently, from \eqref{41} and \eqref{42} we see that
\begin{equation*}
(t-t_{0})^{\frac{1}{3}}\|u_{\varepsilon}(\cdot,t)\|_{L^{3}(\Omega)}\leq9 C_{1}\delta^{\frac{1}{2}}\ \ \ \mbox{for all} \  t\in[t_{0}, t_{0}+2)
\end{equation*}
and
\begin{equation*}
(t-t_{0})^{\frac{1}{6}}\|\nabla v_{\varepsilon}(\cdot,t)\|_{L^{6}(\Omega)}\leq9 C_{4}\delta^{\frac{1}{4}}\ \ \ \mbox{for all} \  t\in[t_{0}, t_{0}+2),
\end{equation*}
which ensure that
\begin{equation}\label{43}
\|u_{\varepsilon}(\cdot,t)\|_{L^{3}(\Omega)}\leq9 C_{1}\delta^{\frac{1}{2}}(t-t_{0})^{-\frac{1}{3}}\leq9 C_{1}\delta^{\frac{1}{2}}\ \ \ \mbox{for all} \  t\in[t_{0}+1, t_{0}+2)
\end{equation}
and
\begin{equation}\label{44}
\|\nabla v_{\varepsilon}(\cdot,t)\|_{L^{6}(\Omega)}\leq  9 C_{4}\delta^{\frac{1}{4}}(t-t_{0})^{-\frac{1}{6}}\leq 9 C_{4}\delta^{\frac{1}{4}}\ \ \ \mbox{for all} \  t\in[t_{0}+1, t_{0}+2).
\end{equation}
Since $C_{1}\delta^{\frac{1}{2}}$ and $C_{4}\delta^{\frac{1}{4}}$ are independent of the choice of $t_0\geq T_{\delta,2}-2:=T_{\delta,1}$, \eqref{43} and \eqref{44} provide a constant $C_{8}=9(C_{1}+C_{4})>0$ fulfilling
\begin{equation*}
\|u_{\varepsilon}(\cdot,t)\|_{L^{3}(\Omega)}\leq C_{8}\delta^{\frac{1}{2}}\ \ \ \mbox{for all} \  t>T_{\delta,2}-1
\end{equation*}
and
\begin{equation*}
\|\nabla v_{\varepsilon}(\cdot,t)\|_{L^{6}(\Omega)}\leq C_{8}\delta^{\frac{1}{4}}\ \ \ \mbox{for all} \  t>T_{\delta,2}-1.
\end{equation*}
This readily leads to \eqref{f1} and \eqref{f2}.
\end{proof}

The improved regularity in Lemma~\ref{f} allows us to derive
eventual bounds for $(u_\varepsilon,v_\varepsilon)$ in
$L^\infty(\Omega)\times W^{1,\infty}(\Omega)$. We also make explicit the dependence of these bounds on the smallness
parameter for the subsequent asymptotic analysis
\begin{lemma}\label{g}
Let $N=4$. There exist $\widetilde{C}>0$ and $\delta_{\star\star}>0$ such that for any $\delta\in(0, \delta_{\star\star}]$, one can find $\lambda_{\delta,3}=\lambda_{\delta,3}(\delta,\chi,\mu,\Omega)\in(0,\frac{1}{16}]$ with the property that whenever $\lambda<\lambda_{\delta,3}$, then there exists $T_{\delta,3}=T_{\delta,3}(\delta,\chi,\mu,\Omega,u_{0},v_{0})>T_{\delta,2}$ so that for any choice of $\varepsilon\in(0,1)$, we have
\begin{equation}\label{g1}
\|u_{\varepsilon}(\cdot,t)\|_{L^{\infty}(\Omega)}\leq \widetilde{C}\delta^{\frac{1}{2}}\ \ \ \mbox{for all} \  t\geq T_{\delta,3}
\end{equation}
and
\begin{equation}\label{g2}
\|\nabla v_{\varepsilon}(\cdot,t)\|_{L^{\infty}(\Omega)}\leq \widetilde{C}\left(\delta^{\frac{1}{2}}+\delta^{\frac{1}{4}}\right)\ \ \ \mbox{for all} \  t\geq T_{\delta,3}.
\end{equation}
Moreover, one has
\begin{equation}\label{g3}
\|v_{\varepsilon}(\cdot,t)\|_{L^{\infty}(\Omega)}\leq \widetilde{C}\left(e^{-(t-T_{\delta,3}+\frac{1}{2})}+\delta^{\frac{1}{2}}\right)\ \ \ \mbox{for all} \  t\geq T_{\delta,3}.
\end{equation}
\end{lemma}
\begin{proof}
For each $\delta\in(0,\delta_\star]$, let $\lambda_{\delta,2}$ be the constant given by Lemma \ref{f}. We define
\begin{equation}\label{45}
\lambda_{\delta,3}:=\min\left\{\lambda_{\delta,2},\frac{1}{16}\right\}
\end{equation}
and assume that $\lambda<\lambda_{\delta,3}$. Then Lemma \ref{f} yields the existence of $T_{\delta,2}>0$ and $C_{1}>0$ such that
\begin{equation}\label{46}
\|u_{\varepsilon}(\cdot,t)\|_{L^{3}(\Omega)}\leq C_{1}\delta^{\frac{1}{2}}\ \ \ \mbox{for all} \  t\geq T_{\delta,2}\ \mbox{and}\ \varepsilon\in(0,1)
\end{equation}
and
\begin{equation}\label{47}
\|\nabla v_{\varepsilon}(\cdot,t)\|_{L^{6}(\Omega)}\leq C_{1}\delta^{\frac{1}{4}}\ \ \ \mbox{for all} \  t\geq T_{\delta,2}\ \mbox{and}\ \varepsilon\in(0,1).
\end{equation}
Let $t_0\geq T_{\delta,3}-2:=T_{\delta,2}$. For any $\varepsilon\in(0,1)$, we also define
\begin{equation}\label{48}
M:=\sup\limits_{\tau\in(t_{0},t_{0}+2]} \|(\tau-t_{0})^{\frac{2}{3}}u_{\varepsilon}(\cdot,\tau)\|_{L^{\infty}(\Omega)}.
\end{equation}
By \eqref{341}, we obtain
\begin{equation}\label{49}
\begin{aligned}
\|\sigma^{\frac{2}{3}}u_{\varepsilon}(\cdot,t_{0}+\sigma)\|_{L^{\infty}(\Omega)}\leq&\underbrace{\sigma^{\frac{2}{3}}\|e^{\sigma\Delta}u_{\varepsilon}(\cdot,t_{0})\|_{L^{\infty}(\Omega)}}\limits_{:=I_{6}}+\underbrace{\lambda_{+}\sigma^{\frac{2}{3}}\int_{0}^{\sigma}\|e^{(\sigma-s)\Delta}u_{\varepsilon}(\cdot,t_{0}+s)\|_{L^{\infty}(\Omega)}ds}\limits_{:=I_{7}}\\
&+\underbrace{\chi\sigma^{\frac{2}{3}}\int_{0}^{\sigma}\|e^{(\sigma-s)\Delta}\nabla \cdot\left(u_{\varepsilon}(\cdot,t_{0}+s)\nabla v_{\varepsilon}(\cdot,t_{0}+s)\right)\|_{L^{\infty}(\Omega)}ds}\limits_{:=I_{8}}\\
\end{aligned}
\end{equation}
for all $\sigma\in(0,2]$ and $\varepsilon\in(0,1)$, where the smoothing estimates for the Neumann heat semigroup $\{e^{t\Delta}\}_{t\geq0}$ and \eqref{46}-\eqref{48} assert that
\begin{equation*}
\begin{aligned}
I_{6}\leq& C_{2}\sigma^{\frac{2}{3}}\left(1+\sigma^{-\frac{2}{3}}\right)\|u_{\varepsilon}(\cdot,t_{0})\|_{L^{3}(\Omega)}\\
\leq& C_{1}C_{2}\left(1+\sigma^{\frac{2}{3}}\right)\delta^{\frac{1}{2}}\leq 3C_{1}C_{2}\delta^{\frac{1}{2}}\ \ \ \mbox{for all} \  \sigma\in(0,2]\ \mbox{and}\ \varepsilon\in(0,1)
\end{aligned}
\end{equation*}
and
\begin{equation*}
\begin{aligned}
I_{7}\leq& \lambda_{+}\sigma^{\frac{2}{3}}\int_{0}^{\sigma}\|u_{\varepsilon}(\cdot,t_{0}+s)\|_{L^{\infty}(\Omega)}ds\\
\leq& \lambda_{+}\sigma^{\frac{2}{3}}M\int_{0}^{\sigma}s^{-\frac{2}{3}}ds\leq 6\lambda_{+}M\ \ \ \mbox{for all} \  \sigma\in(0,2]\ \mbox{and}\ \varepsilon\in(0,1)
\end{aligned}
\end{equation*}
as well as
\begin{equation*}
\begin{aligned}
I_{8}\leq& \chi C_{3}\sigma^{\frac{2}{3}}\int_{0}^{\sigma}\left(1+(\sigma-s)^{-\frac{5}{6}}\right)\|u_{\varepsilon}(\cdot,t_{0}+s)\nabla v_{\varepsilon}(\cdot,t_{0}+s)\|_{L^{6}(\Omega)}ds\\
\leq& \chi C_{3}\sigma^{\frac{2}{3}}\int_{0}^{\sigma}\left(1+(\sigma-s)^{-\frac{5}{6}}\right)\|u_{\varepsilon}(\cdot,t_{0}+s)\|_{L^{\infty}(\Omega)}\|\nabla v_{\varepsilon}(\cdot,t_{0}+s)\|_{L^{6}(\Omega)}ds\\
\leq& \chi C_{1}C_{3}\delta^{\frac{1}{4}}M\cdot\left\{\sigma^{\frac{2}{3}}\int_{0}^{\sigma}s^{-\frac{2}{3}}\left(1+(\sigma-s)^{-\frac{5}{6}}\right) ds\right\}\\
\leq& \chi C_{1}C_{3}C_{4}\delta^{\frac{1}{4}}M\ \ \ \mbox{for all} \  \sigma\in(0,2]\ \mbox{and}\ \varepsilon\in(0,1).
\end{aligned}
\end{equation*}
with some constants $C_{i}>0$ ($i=2, 3, 4$) thanks to $\sigma^{\frac{2}{3}}\int_{0}^{\sigma}s^{-\frac{2}{3}}\left(1+(\sigma-s)^{-\frac{5}{6}}\right)ds$ is finite.
In consequence, \eqref{49} entails that
\begin{equation}\label{50}
M\leq3C_{1}C_{2}\delta^{\frac{1}{2}}+6\lambda_{+}M+\chi C_{1}C_{3}C_{4}\delta^{\frac{1}{4}}M\ \ \ \mbox{for all} \  \sigma\in(0,2]\ \mbox{and}\ \varepsilon\in(0,1).
\end{equation}
Now we define
\begin{equation}\label{51}
\delta_{\star\star}:=\min\left\{\delta_{\star}, \frac{1}{(4\chi C_{1}C_{3}C_{4})^{4}}\right\}
\end{equation}
and fix any $\delta\in(0,\delta_{\star\star}]$. Since \eqref{45} and \eqref{51} ensure $6\lambda_{+}M\leq\frac{M}{2}$ and $\chi C_{1}C_{3}C_{4}\delta^{\frac{1}{4}}M\leq\frac{M}{4}$, we deduce from \eqref{50} that
\begin{equation*}
(t-t_{0})^{\frac{2}{3}}\|u_{\varepsilon}(\cdot,t)\|_{L^{\infty}(\Omega)}\leq12 C_{1}C_{2}\delta^{\frac{1}{2}}\ \ \ \mbox{for all} \  t\in[t_{0}, t_{0}+2]
\end{equation*}
which implies
\begin{equation}\label{52}
\|u_{\varepsilon}(\cdot,t)\|_{L^{\infty}(\Omega)}\leq12 C_{1}C_{2}\delta^{\frac{1}{2}}(t-t_{0})^{-\frac{2}{3}}\leq12 C_{1}C_{2}\delta^{\frac{1}{2}}\ \ \ \mbox{for all} \  t\in[t_{0}+1, t_{0}+2].
\end{equation}
Observing that $12 C_{1}C_{2}\delta^{\frac{1}{2}}$ is independent of $t_0$ for all $t_0\geq T_{\delta,3}-2:=T_{\delta,2}$, from \eqref{52} we readily obtain
\begin{equation}\label{53}
\|u_{\varepsilon}(\cdot,t)\|_{L^{\infty}(\Omega)}\leq C_{5}\delta^{\frac{1}{2}}\ \ \ \mbox{for all} \  t\geq T_{\delta,3}-1\ \mbox{and}\ \varepsilon\in(0,1)
\end{equation}
with $C_{5}>0$. Moreover, let $t_{0}=T_{\delta,3}-1$ and it follows from the Duhamel representation that for all $t\geq t_{0}+1$,
\begin{equation*}
v_{\varepsilon}(\cdot,t)= e^{(t-t_{0})(\Delta-1)}v_{\varepsilon}(\cdot,t_{0})+\int_{t_{0}}^{t} e^{(t-s)(\Delta-1)}u_{\varepsilon}(\cdot,s)ds,
\end{equation*}
which allows us to derive
\begin{equation}\label{54}
\|\nabla v_{\varepsilon}(\cdot,t)\|_{L^{\infty}(\Omega)}
\leq\underbrace{\|\nabla e^{(t-t_{0})(\Delta-1)} v_{\varepsilon}(\cdot,t_{0})\|_{L^{\infty}(\Omega)}}\limits_{:=I_{9}}+\underbrace{\int_{t_{0}}^{t}\|\nabla e^{(t-s)(\Delta-1)}u_{\varepsilon}(\cdot,s)\|_{L^{\infty}(\Omega)}ds}\limits_{:=I_{10}}
\end{equation}
and
\begin{equation}\label{55}
\| v_{\varepsilon}(\cdot,t)\|_{L^{\infty}(\Omega)}
\leq\underbrace{\|e^{(t-t_{0})(\Delta-1)} v_{\varepsilon}(\cdot,t_{0})\|_{L^{\infty}(\Omega)}}\limits_{:=I_{11}}+\underbrace{\int_{t_{0}}^{t}\| e^{(t-s)(\Delta-1)}u_{\varepsilon}(\cdot,s)\|_{L^{\infty}(\Omega)}ds}\limits_{:=I_{12}}
\end{equation}
for all $t\geq t_{0}+1$ and $\varepsilon\in(0,1)$. Here, Lemma \ref{p2} together with \eqref{d4}, \eqref{47} and \eqref{53} ensures that
\begin{equation}\label{56}
\begin{aligned}
I_{9}\leq& C_{6}\left(1+(t-t_{0})^{-\frac{1}{3}}\right)\|\nabla v_{\varepsilon}(\cdot,t_{0})\|_{L^{6}(\Omega)}\\
\leq& 2C_{1}C_{6}\delta^{\frac{1}{4}}\ \ \ \mbox{for all} \  t\geq t_{0}+1\ \mbox{and}\ \varepsilon\in(0,1)
\end{aligned}
\end{equation}
and
\begin{equation}\label{57}
\begin{aligned}
I_{10}\leq& C_{7}\int_{t_{0}}^{t}\left(1+(t-s)^{-\frac{1}{2}}\right)e^{-(t-s)}\| u_{\varepsilon}(\cdot,s)\|_{L^{\infty}(\Omega)}ds\\
\leq& C_{5}C_{7}C_{8}\delta^{\frac{1}{2}}\ \ \ \mbox{for all} \  t\geq t_{0}+1\ \mbox{and}\ \varepsilon\in(0,1)
\end{aligned}
\end{equation}
as well as
\begin{equation}\label{58}
\begin{aligned}
I_{11}\leq& e^{-(t-t_{0})}C_{9}\left(1+(t-t_{0})^{-1}\right)\|v_{\varepsilon}(\cdot,t_{0})\|_{L^{2}(\Omega)}\\
\leq& 2e^{-(t-t_{0})}C_{9}\sqrt{M_{0}}\ \ \ \mbox{for all} \  t\geq t_{0}+1\ \mbox{and}\ \varepsilon\in(0,1)
\end{aligned}
\end{equation}
and
\begin{equation}\label{59}
\begin{aligned}
I_{12}\leq& C_{10}\int_{t_{0}}^{t}\left(1+(t-s)^{-\frac{1}{2}}\right)e^{-(t-s)}\| u_{\varepsilon}(\cdot,s)\|_{L^{\infty}(\Omega)}ds\\
\leq& C_{5}C_{8}C_{10}\delta^{\frac{1}{2}}\ \ \ \mbox{for all} \  t\geq t_{0}+1\ \mbox{and}\ \varepsilon\in(0,1).
\end{aligned}
\end{equation}
with $C_{i}>0$ ($i=6, ..., 10$). Substituting \eqref{56}-\eqref{57} into \eqref{54} and \eqref{58}-\eqref{59} into \eqref{55}, we see that for all $t\geq T_{\delta,3}$ and $\varepsilon\in(0,1)$,
\begin{equation}\label{60}
\|\nabla v_{\varepsilon}(\cdot,t)\|_{L^{\infty}(\Omega)}
\leq2C_{1}C_{6}\delta^{\frac{1}{4}}+C_{5}C_{7}C_{8}\delta^{\frac{1}{2}}
\end{equation}
and
\begin{equation}\label{61}
\| v_{\varepsilon}(\cdot,t)\|_{L^{\infty}(\Omega)}
\leq\frac{2}{e}\cdot e^{-(t-T_{\delta,3})}C_{9}\sqrt{M_{0}}+C_{5}C_{8}C_{10}\delta^{\frac{1}{2}}.
\end{equation}
Therefore, the claim \eqref{g1}-\eqref{g3} results from \eqref{53} and \eqref{60}-\eqref{61} upon choosing $\widetilde{C}=\max\{C_{5},2C_{1}C_{6},C_{5}C_{7}C_{8},C_{9}\sqrt{M_{0}},C_{5}C_{8}C_{10}\}+1$.
\end{proof}

In particular, the eventual $L^\infty$-bound for
$u_\varepsilon$ can be made arbitrarily small by taking
the growth rate sufficiently small.
\begin{corollary}\label{h}
For every $\nu>0$ there exists
$\lambda_\nu=\lambda_\nu(\nu,\chi,\mu,\Omega)>0$ with the property that if $\lambda<\lambda_\nu$, then one can pick $T_{\nu}=T_{\nu}(\nu,\chi,\mu,\Omega,u_{0},v_{0})>0$ such that
\begin{equation}\label{h1}
\|u_\varepsilon(\cdot,t)\|_{L^\infty(\Omega)}\le \nu
\ \ \ \mbox{for all} \  t\geq T_{\nu}\ \mbox{and}\ \varepsilon\in(0,1).
\end{equation}
\end{corollary}
\begin{proof}
Since only \eqref{g1} is needed here, we may replace
$\widetilde{C}$ with $C_5$ in this estimate,
as seen from \eqref{53}. For any fixed $\nu>0$, we apply Lemma \ref{g} with
\begin{equation*}
\delta=\min\left\{\left(\frac{\nu}{2C_{5}}\right)^{2}, \frac{\delta_{\star\star}}{2}\right\},
\end{equation*}
Let $\lambda_\nu(\nu,\chi,\mu,\Omega):=\lambda_{\delta,3}$. Then whenever $\lambda<\lambda_\nu$, we can find $T_{\nu}(\nu,\chi,\mu,\Omega,u_{0},v_{0})=T_{\delta,3}+1>0$ such that
\begin{equation*}
\|u_{\varepsilon}(\cdot,t)\|_{L^{\infty}(\Omega)}\leq C_{5}\delta^{\frac{1}{2}}<\nu\ \ \ \mbox{for all} \  t\geq T_{\nu}\ \mbox{and}\ \varepsilon\in(0,1).
\end{equation*}
This is our desired estimate \eqref{h1}.
\end{proof}

With the bounds obtained in Lemma~\ref{g}, standard
parabolic H\"older estimates and Schauder theory
lead to the higher regularity stated below. The argument is analogous to those in \cite[Lemma~2.12]{D2026} and \cite[Lemma~8.1]{W2022} and we omit the details.
\begin{lemma}\label{i}
Under the conditions in Lemma \ref{g}, there exists $\widetilde{T}_{1}>0$ and $\alpha_{1}\in(0,1)$ such that
\begin{equation}\label{i1}
\|u_{\varepsilon}\|_{C^{2+\alpha_{1},1+\frac{\alpha_{1}}{2}}(\overline{\Omega}\times[t,t+1])}\leq C\ \ \ \mbox{for all}\ \ t>\widetilde{T}_{1}
\end{equation}
and
\begin{equation}\label{i2}
\|v_{\varepsilon}\|_{C^{2+\alpha_{1},1+\frac{\alpha_{1}}{2}}(\overline{\Omega}\times[t,t+1])}\leq C\ \ \ \mbox{for all}\ \ t>\widetilde{T}_{1}
\end{equation}
with a constant $C>0$.
\end{lemma}

The following convergence result follows from the
Arzel\`a--Ascoli theorem, as in \cite[Lemma~8.2]{W2022}.
We omit the standard compactness argument.
\begin{lemma}\label{ii}
Let the assumptions of Lemma \ref{i} hold, and Let $(u_\varepsilon,v_\varepsilon)$ be a global solution to \eqref{model111}. Then there exist $\alpha_{2}\in(0,\alpha_{1})$ and $\widetilde{T}_{2}>\widetilde{T}_{1}$ with the property that given any
$(\varepsilon_j)_{j\in\mathbb{N}}\subset(0,1)$ satisfying $\varepsilon_j\to0$ as $j\to\infty$, one can find
\begin{equation}\label{iii1}
\left\{\aligned
&u\in C^{2+\alpha_{2},1+\frac{\alpha_{2}}{2}}(\overline{\Omega}\times[t,t+1]),
&&\\
&v\in C^{2+\alpha_{2},1+\frac{\alpha_{2}}{2}}(\overline{\Omega}\times[t,t+1]),
\endaligned\right.
\end{equation}
and a subsequence $(\varepsilon_{j_k})_{k\in\mathbb{N}}$ such that
\begin{equation}\label{iii2}
\left\{\aligned
&u_{\varepsilon}\rightarrow u\ \ \mbox{in}\ C_{loc}^{2,1}(\overline{\Omega}\times[\widetilde{T}_{2},\infty)),
&&\\
&v_{\varepsilon}\rightarrow v\ \ \mbox{in}\ C_{loc}^{2,1}(\overline{\Omega}\times[\widetilde{T}_{2},\infty)),
\endaligned\right.
\end{equation}
as $\varepsilon=\varepsilon_{jk}\searrow0$.
\end{lemma}

The smallness estimate in Corollary~\ref{h} is also preserved
in the limit.
\begin{corollary}\label{iii}
Let $(u,v)$ be the global weak solution constructed in Proposition \ref{prop1}, and assume that the hypotheses of Theorem \ref{theorem1} hold. For every $\nu>0$, let $\lambda_\nu>0$ be given by Corollary \ref{h}. If $\lambda<\min\{\lambda_0,\lambda_\nu\}$, then one can find $\widehat{T}_{\nu}>T_0$ such that
\begin{equation}\label{iiii1}
\|u(\cdot,t)\|_{L^\infty(\Omega)}\leq\nu
\ \ \ \mbox{for all}\ \ t\geq\widehat{T}_{\nu}.
\end{equation}
\end{corollary}
\begin{proof}
Fix $\nu>0$ and suppose that $\lambda<\min\{\lambda_0,\lambda_\nu\}$. According to Corollary \ref{h}, there exists $T_{\nu}>0$ independent of $\varepsilon$ such that
\begin{equation}\label{iiii2}
\|u_{\varepsilon}(\cdot,t)\|_{L^\infty(\Omega)}\leq\nu
\ \ \ \mbox{for all}\ \ t\geq T_{\nu}\ \mbox{and}\ \varepsilon\in(0,1).
\end{equation}
Let $\widehat{T}_{\nu}:=\max\{T_{\nu},T_0\}+1$. By the almost everywhere convergence in Proposition \ref{prop1} and the eventual regularity asserted in Theorem \ref{theorem1}, passing to the limit in \eqref{iiii2} yields \eqref{iiii1}.
\end{proof}
\par\noindent\textbf{Proof of Theorem~\ref{theorem1}.}
Fix $\delta=\frac{\delta_{\star\star}}{2}$ and let $\lambda_0:=\lambda_{\delta,3}>0$ be as in Lemma~\ref{g}. Combining Lemma \ref{ii} with Proposition \ref{prop1}, we conclude that there exists $T_{0}:=\widetilde{T}_{2}+1>0$ such that the global weak solution $(u,v)$ established in Proposition \ref{prop1} enjoys the eventual smoothness properties stated in Theorem \ref{theorem1} and indeed solves \eqref{model11} classically in $\overline{\Omega}\times[T_{0},\infty)$. Finally, the uniform bounds derived in Lemma \ref{i} are preserved in the limit, thereby establishing \eqref{tt2} for any fixed $\alpha\in(0,\alpha_{1})$.

\section{Sharp asymptotics}
\indent

In this section, we shall prove Theorems~\ref{theorem2}--\ref{theorem4}
on the sharp asymptotic behaviour of global weak solutions. Throughout this section, we assume that the conditions in Theorem \ref{theorem1} hold. Let $(u,v)$ denote the global weak solution constructed in Proposition \ref{prop1}, and let $\lambda_{0}$ and $T_{0}$ be as in Theorem \ref{theorem1}.

We first establish strict positivity of the cell density for all $t>T_0$. This will justify taking
logarithms and reciprocals of its spatial average.

\setcounter{equation}{0}
\begin{lemma}\label{j}
Suppose that $\lambda<\lambda_{0}$. Then
\begin{equation}\label{j1}
u(x,t)>0\ \ \ \mbox{in} \ \overline{\Omega}\times(T_{0},\infty).
\end{equation}
\end{lemma}

\begin{proof}
According to proposition \ref{prop1}, the strong maximum principle ensures that $u_{\varepsilon}$ is strictly positive in $\overline{\Omega}\times(0,\infty)$. Then testing the first equation in \eqref{model111} against $\frac{1}{u_\varepsilon}$ gives

\begin{equation}\label{62}
\begin{aligned}
\frac{d}{dt}\int_\Omega\ln u_\varepsilon
={}&\int_\Omega\frac{|\nabla u_\varepsilon|^2}{u_\varepsilon^2}-\chi\int_\Omega\frac{\nabla u_\varepsilon}{u_\varepsilon}\cdot\nabla v_\varepsilon+\lambda|\Omega|-\mu\int_\Omega u_\varepsilon
-\varepsilon\int_\Omega u_\varepsilon^{\theta-1}\\
\ge{}&-\frac{\chi^2}{4}\int_\Omega|\nabla v_\varepsilon|^2
 +\lambda|\Omega|-\mu\int_\Omega u_\varepsilon
 -\varepsilon\int_\Omega u_\varepsilon^{\theta-1}\ \ \ \mbox{for all} \  t>0\ \mbox{and}\ \varepsilon\in(0,1).
\end{aligned}
\end{equation}
Multiplying the second equation in \eqref{model111} by $v_{\varepsilon}$, an application of Young's inequality shows
\begin{equation}\label{63}
\frac{1}{2}\frac{\rm d}{{\rm d}t}\int_{\Omega}v_{\varepsilon}^{2}+\int_{\Omega}|\nabla v_{\varepsilon}|^{2}=-\int_{\Omega}v_{\varepsilon}^{2}+\int_{\Omega}u_{\varepsilon}v_{\varepsilon}
\leq-\frac{1}{2}\int_{\Omega}v_{\varepsilon}^{2}+\frac{1}{2}\int_{\Omega}u_{\varepsilon}^{2}
\end{equation}
for all $t>0$ and $\varepsilon\in(0,1)$. The combination of \eqref{62} and \eqref{63} yields
\begin{equation}\label{64}
\begin{aligned}
\frac{d}{dt}\left\{-\int_{\Omega}\ln u_\varepsilon+\frac{\chi^{2}}{8}\int_{\Omega}v_{\varepsilon}^{2}\right\}
\leq&-\lambda|\Omega|+\mu\int_\Omega u_\varepsilon+\varepsilon\int_\Omega u_\varepsilon^{\theta-1}+\frac{\chi^{2}}{8}\int_{\Omega}u_{\varepsilon}^{2}\\
\leq&\left(\mu+\frac{\chi^{2}}{8}\right)\int_{\Omega}u_{\varepsilon}^{2}+\varepsilon\int_\Omega u_\varepsilon^{\theta}+(\mu+1+\lambda_{-})|\Omega|
\end{aligned}
\end{equation}
for all $t>0$ and $\varepsilon\in(0,1)$, because $\varepsilon u_\varepsilon^{\theta-1}\le1+\varepsilon u_\varepsilon^\theta$. Based on \eqref{d0}, we integrate \eqref{6} in time to derive
\begin{equation*}
\int_0^t\left(\mu \int_{\Omega} u_{\varepsilon}^{2}+\varepsilon \int_{\Omega} u_\varepsilon^{\theta}\right)\leq \lambda_{+}\max\left\{\int_{\Omega}u_{0},\frac{\lambda_{+}|\Omega|}{\mu}\right\}\cdot t+\int_{\Omega} u_0\ \ \ \mbox{for all} \  t>0\ \mbox{and}\ \varepsilon\in(0,1).
\end{equation*}
Accordingly, a direct integration of \eqref{64} implies that for all $T>T_{0}$,
\begin{equation}\label{65}
\begin{aligned}
&-\int_{\Omega}\ln u_{\varepsilon}(\cdot,t)+\frac{\chi^{2}}{8}\int_{\Omega}v_{\varepsilon}^{2}(\cdot,t)\\
\leq&-\int_{\Omega}\ln u_0+\frac{\chi^{2}}{8}\int_{\Omega}v_0^2
+\int_0^t\left\{\left(\mu+\frac{\chi^{2}}{8}\right)\int_{\Omega}u_{\varepsilon}^{2}+\varepsilon\int_\Omega u_\varepsilon^{\theta}\right\}
+(\mu+1+\lambda_{-})|\Omega|\cdot t\\
\leq& C_{1}(T)\ \ \ \mbox{for all} \  t\in(0,T)\ \mbox{and}\ \varepsilon\in(0,1)
\end{aligned}
\end{equation}
with some $C_{1}(T)>0$. Since the Jensen's inequality warrants that
\begin{equation*}
\frac{1}{|\Omega|}
\int_{\Omega}\ln u_\varepsilon(\cdot,t)
\leq
\ln\left(
\frac{1}{|\Omega|}
\int_{\Omega}u_\varepsilon(\cdot,t)\right)\ \ \ \mbox{for all} \  t>0\ \mbox{and}\ \varepsilon\in(0,1),
\end{equation*}
from \eqref{65} we infer that
\begin{equation*}
\begin{aligned}
\int_{\Omega}u_\varepsilon(\cdot,t)
\geq{}&
|\Omega|
\exp\left(
\frac{1}{|\Omega|}
\int_{\Omega}\ln u_\varepsilon(\cdot,t)
\right)\\
\geq{}&
|\Omega|
\exp\left(
-\frac{C_1(T)}{|\Omega|}
\right)>0
\end{aligned}
\end{equation*}
for all $t\in(0,T]$ and $\varepsilon\in(0,1)$.
In particular, since $T>T_0$, we have
\begin{equation*}
\int_{\Omega}u_\varepsilon(\cdot,T_0)
\geq
|\Omega|\exp\left(
-\frac{C_1(T)}{|\Omega|}
\right)>0
\ \ \ \mbox{for all}\ \
\varepsilon\in(0,1).
\end{equation*}
Applying Lemma \ref{ii} to the sequence selected in Proposition \ref{prop1}, we further have
\begin{equation*}
\int_{\Omega} u(\cdot,T_{0})>0,
\end{equation*}
so that
\begin{equation}\label{67}
u(\cdot,T_{0})\geq, \not\equiv0.
\end{equation}
In view of \eqref{67} and the eventual regularity asserted in Theorem \ref{theorem1}, another application of the strong maximum principle immediately yields \eqref{j1}.
\end{proof}

The following lemma gives the evolution equations for the spatial averages.
\begin{lemma}\label{k}
The spatial average $\bar u$ satisfies
\begin{equation}\label{k1}
\frac{d}{dt}\bar u(t)=\lambda\bar u(t)
-\frac{\mu}{|\Omega|}\int_\Omega u^2
\ \ \ \mbox{for all}\ t\geq T_{0},
\end{equation}
or equivalently,
\begin{equation}\label{k2}
\frac{d}{dt}\bar u(t)
=\lambda\bar u(t)-\mu\bar u^2(t)
-\frac{\mu}{|\Omega|}\int_\Omega(u-\bar u)^2\ \ \ \mbox{for all}\ t\geq T_{0}.
\end{equation}
Moreover, the spatial average $\bar v$ satisfies
\begin{equation}\label{k3}
\frac{d}{dt}\bar v(t)=-\bar v(t)+\bar u(t).
\end{equation}
\end{lemma}
\begin{proof}
Integrating the first equation in \eqref{model11}, we obtain
\begin{equation*}
\begin{aligned}
\frac{d}{dt}\bar u(t)=&\lambda\bar u(t)
-\frac{\mu}{|\Omega|}\int_\Omega u^2\\
=&\lambda\bar u(t)-\frac{\mu}{|\Omega|}\left(\int_\Omega (u-\bar u)^2+|\Omega|\bar u^2\right)\ \ \ \mbox{for all}\ t\geq T_{0},
\end{aligned}
\end{equation*}
and thereby establishes \eqref{k1} and \eqref{k2}. Moreover, integrating the second equation in \eqref{model11} readily yields \eqref{k3}.
\end{proof}

We proceed to show the exponential decay of the deviations from the spatial averages in
$L^\infty(\Omega)\times W^{1,\infty}(\Omega)$, which is essential for the subsequent analysis.
\begin{lemma}\label{l}
There exists $\lambda_\star\in(0,\lambda_0]$ with the
property that whenever
$\lambda<\lambda_\star$, one can find $T_{1}>T_0$ and $C>0$ such that
\begin{equation}\label{l1}
\|u(\cdot,t)-\bar u(t)\|_{L^\infty(\Omega)}
+\|v(\cdot,t)-\bar v(t)\|_{W^{1,\infty}(\Omega)}\leq Ce^{-\frac{\rho_1}{12}t}
\ \ \ \mbox{for all}\ t\geq T_{1}.
\end{equation}
\end{lemma}

\begin{proof}
According to \eqref{tt2}, we have
\begin{equation}\label{73}
\| u(\cdot,t)\|_{W^{1,\infty}(\Omega)}
+\|v(\cdot,t)\|_{W^{1,\infty}(\Omega)}\leq C_1
\ \ \ \mbox{for all}\ t\geq T_1
\end{equation}
with $C_1>0$. Choosing $\nu=\frac{\sqrt{\rho_1}}{2\chi}$, we define
\begin{equation*}
\lambda_{\star}:=\min\left\{\lambda_{0},\lambda_{\nu},\frac{\rho_{1}}{8}\right\},
\end{equation*}
where $\lambda_{\nu}$ is defined in Corollary \ref{iii}. Suppose that $\lambda<\lambda_\star$. Then Corollary \ref{iii} guarantees the existence of $T'_{1}>T_0$ such that
\begin{equation}\label{22222}
\|u(\cdot,t)\|_{L^\infty(\Omega)}\leq\nu
\ \ \ \mbox{for all}\ \ t\geq T'_{1}.
\end{equation}
Since $\int_\Omega(u-\bar{u})=0$ and $\nabla(u-\bar{u})=\nabla u$,
from the first equation in \eqref{model11} we see that
\begin{equation}\label{68}
\begin{aligned}
\frac12\frac{\rm d}{{\rm d}t}\int_\Omega(u-\bar{u})^2
=&\int_\Omega(u-\bar{u})(u_t-\frac{\partial}{\partial t}\bar{u})\\
=&\int_\Omega(u-\bar{u})
\bigl(\Delta u-\chi\nabla\cdot(u\nabla v)+\lambda u-\mu u^2\bigr)\\
=&-\int_\Omega|\nabla u|^2
+\chi\int_\Omega u\nabla u\cdot\nabla v
+\lambda\int_\Omega(u-\bar{u})^2
-\mu\int_\Omega(u-\bar{u})^2(u+\bar{u})
\end{aligned}
\end{equation}
for all $t\geq T'_{1}$. Here we may apply Young's inequality and \eqref{22222} to derive
\begin{equation*}
\chi\int_\Omega u\nabla u\cdot\nabla v
\leq\frac12\int_\Omega|\nabla u|^2
+\frac{\chi^2\nu^2}{2}\int_\Omega|\nabla v|^2
\ \ \ \mbox{for all}\ t\geq T'_{1}
\end{equation*}
with $\nu=\frac{\sqrt{\rho_1}}{2\chi}$, which updates \eqref{68} as
\begin{equation}\label{69}
\frac12\frac{\rm d}{{\rm d}t}\int_\Omega(u-\bar{u})^2\\
\leq-\frac12\int_\Omega|\nabla u|^2+\frac{\chi^2\nu^2}{2}\int_\Omega|\nabla v|^2+\lambda\int_\Omega(u-\bar{u})^2
\ \ \ \mbox{for all}\ t\geq T'_{1}.
\end{equation}
Likewise, the second equation in \eqref{model11} and Young's inequality allow us to obtain
\begin{equation}\label{70}
\begin{aligned}
\frac12\frac{\rm d}{{\rm d}t}\int_\Omega(v-\bar{v})^2=&-\int_\Omega|\nabla v|^2-\int_\Omega(v-\bar{v})^2
+\int_\Omega(u-\bar{u})(v-\bar{v})\\
\leq&-\int_\Omega|\nabla v|^2
-\frac12\int_\Omega(v-\bar{v})^2
+\frac12\int_\Omega(u-\bar{u})^2
\ \ \ \mbox{for all}\ t\geq T'_{1}.
\end{aligned}
\end{equation}
The combination of \eqref{69} and \eqref{70} shows
\begin{equation}\label{71}
\begin{aligned}
&\frac12\frac{\rm d}{{\rm d}t}
\left\{\int_\Omega(u-\bar{u})^2
+\frac{\rho_1}{4}\int_\Omega(v-\bar{v})^2\right\}\\
\leq&\underbrace{-\frac12\int_\Omega|\nabla u|^2
+\left(\lambda+\frac{\rho_1}{8}\right)\int_\Omega(u-\bar{u})^2}_{:=J_1}+\underbrace{\left(\frac{\chi^2\nu^2}{2}-\frac{\rho_1}{4}\right)
\int_\Omega|\nabla v|^2-\frac{\rho_1}{8}\int_\Omega(v-\bar{v})^2}_{:=J_2}
\end{aligned}
\end{equation}
for all $t\geq T'_{1}$, where by Lemma \ref{p1} the selections of $\nu$ and $\lambda$,
\begin{equation*}
\begin{aligned}
J_1\leq&\left(-\frac{\rho_1}{2}+\lambda+\frac{\rho_1}{8}\right)
\int_\Omega(u-\bar{u})^2\\
\leq&-\frac{\rho_1}{4}\int_\Omega(u-\bar{u})^2
\ \ \ \mbox{for all}\ t\geq T'_{1}
\end{aligned}
\end{equation*}
and
\begin{equation*}
\begin{aligned}
J_2(t)=&-\frac{\rho_1}{8}\int_\Omega|\nabla v|^2
-\frac{\rho_1}{8}\int_\Omega(v-\bar{v})^2\\
\leq&-\left(\frac{\rho_1^2}{8}+\frac{\rho_1}{8}\right)
\int_\Omega(v-\bar{v})^2\\
\leq&-\frac{\rho_1^2}{16}\int_\Omega(v-\bar{v})^2
\ \ \ \mbox{for all}\ t\geq T'_{1}.
\end{aligned}
\end{equation*}
Therefore, \eqref{71} implies
\begin{equation}\label{7777}
\frac{\rm d}{{\rm d}t}
\left\{\int_\Omega(u-\bar{u})^2+\frac{\rho_1}{4}\int_\Omega(v-\bar{v})^2\right\}+\frac{\rho_1}{2}\left\{\int_\Omega(u-\bar{u})^2+\frac{\rho_1}{4}\int_\Omega(v-\bar{v})^2\right\}\leq0
\end{equation}
for all $t\geq T'_{1}$. The straightforward time integration of \eqref{7777} gives rise to
\begin{equation}\label{77}
\begin{aligned}
&\int_\Omega(u-\bar{u})^2+\frac{\rho_1}{4}\int_\Omega(v-\bar{v})^2\\
\leq& \left\{\int_\Omega\bigl(u(\cdot,T'_1)-\bar{u}(T'_1)\bigr)^2+\frac{\rho_1}{4}\int_\Omega\bigl(v(\cdot,T'_1)-\bar{v}(T'_1)\bigr)^2\right\}\cdot e^{-\frac{\rho_1}{2}(t-T'_{1})}\ \ \ \mbox{for all}\ t\geq T'_{1},
\end{aligned}
\end{equation}
which yields the existence of $C_2>0$ such that
\begin{equation}\label{72}
\|u(\cdot,t)-\bar{u}(t)\|_{L^2(\Omega)}
+\|v(\cdot,t)-\bar{v}(t)\|_{L^2(\Omega)}
\leq C_2e^{-\frac{\rho_1}{4}t}
\ \ \ \mbox{for all}\ t\geq T'_{1}.
\end{equation}
thanks to \eqref{73}. Then the Gagliardo-Nirenberg inequality, \eqref{72} and \eqref{73} provide $C_i>0$ $(i=3, 4, 5)$ such that
\begin{equation}\label{74}
\begin{aligned}
\|u(\cdot,t)-\bar{u}(t)\|_{L^\infty(\Omega)}
\leq &C_3\|\nabla u(\cdot,t)\|_{L^\infty(\Omega)}^{\frac23}
\|u(\cdot,t)-\bar{u}(t)\|_{L^2(\Omega)}^{\frac13}
+C_3\|u(\cdot,t)-\bar{u}(t)\|_{L^2(\Omega)}\\
\leq&C_4\|u(\cdot,t)-\bar{u}(t)\|_{L^2(\Omega)}^{\frac13}
+C_3\|u(\cdot,t)-\bar{u}(t)\|_{L^2(\Omega)}\\
\leq&C_5e^{-\frac{\rho_1}{12}t}
\ \ \ \mbox{for all}\ t\geq T'_{1}.
\end{aligned}
\end{equation}
The same argument applied to $v-\bar v$ also gives
\begin{equation}\label{75}
\|v(\cdot,t)-\bar{v}(t)\|_{L^\infty(\Omega)}
\leq C_6e^{-\frac{\rho_1}{12}t}
\ \ \ \mbox{for all}\ t\geq T'_{1}
\end{equation}
with $C_6>0$. In view of a variation-of-constants representation associated with the second equation in \eqref{model11}, we obtain that
\begin{equation*}
\begin{aligned}
&v(\cdot,t)-\bar{v}(t)\\
= & e^{(t-T'_1)(\Delta-1)}\bigl(v(\cdot,T'_1)-\bar{v}(T'_1)\bigr)+\int_{T'_1}^te^{(t-s)(\Delta-1)}\bigl(u(\cdot,s)-\bar{u}(s)\bigr)ds\ \ \ \mbox{for all}\ t\geq T'_{1}.
\end{aligned}
\end{equation*}
The smoothing estimates for the Neumann heat semigroup together with \eqref{74}-\eqref{75} therefore ensures that
\begin{equation}\label{777}
\begin{aligned}
\|\nabla v(\cdot,t)\|_{L^\infty(\Omega)}\leq&C_{7}\bigl(1+(t-T'_1)^{-\frac12}\bigr)
e^{-(\rho_1+1)(t-T'_1)}\|v(\cdot,T'_1)-\bar{v}(T'_1)\|_{L^{\infty}(\Omega)}
\\&+C_8\int_{T'_1}^t\left(1+(t-s)^{-\frac12}\right)
e^{-\rho_1(t-s)}\|u(\cdot,s)-\bar{u}(s)\|_{L^{\infty}(\Omega)}\,ds\\
\leq&2C_6C_7e^{-\frac{\rho_1}{12}t}
+C_{5}C_8e^{-\frac{\rho_1}{12}t}\int_{T'_1}^t
\left(1+(t-s)^{-\frac12}\right)e^{-\frac{11\rho_1}{12}(t-s)}\,ds\\
\leq&C_9e^{-\frac{\rho_1}{12}t}
\ \ \ \mbox{for all}\ t\geq T'_{1}+1
\end{aligned}
\end{equation}
with some $C_i>0$ $(i=7, 8,9)$. With $T_{1}=T'_{1}+1$, \eqref{l1} follows from \eqref{74}, \eqref{75} and \eqref{777}.
\end{proof}

For zero and small positive growth, the equation for $\frac{1}{\bar u}$ requires control of the variance relative to $\bar u^2$. The next theorem provides the required estimate.
\begin{lemma}\label{m}
Let $\lambda_\star$ and $T_{1}$ be given by Lemma \ref{l}. Then there exists $\lambda_{\star\star}\in(0,\lambda_\star]$ with the
property that whenever $0\leq\lambda<\lambda_{\star\star}$, one can find $T_{2}>T_1$ and $C>0$ such that
\begin{equation}\label{m1}
\frac{\frac1{|\Omega|}\int_\Omega(u-\bar u)^2}{\bar u(t)^2}
\leq Ce^{-\frac{3\rho_1}{8}t}
\ \ \ \mbox{for all}\ t\geq T_{2}.
\end{equation}
\end{lemma}
\begin{proof}
Fixing $\tilde{\nu}>0$ such that
\begin{equation}\label{1144}
\tilde{\nu}\leq\min\left\{\frac{\sqrt{\rho_1}}{2\chi},
\frac{\rho_1}{16\mu}\right\},
\end{equation}
we define
\begin{equation*}
\lambda_{\star\star}:=\min\left\{\lambda_{\star},\lambda_{\tilde{\nu}}\right\},
\end{equation*}
where $\lambda_{\tilde{\nu}}$ is specified in Corollary \ref{iii}. Let $0\leq\lambda<\lambda_{\star\star}$. By \eqref{k1} and Corollary \ref{iii}, there exists $T_2>T_1$ such that
\begin{equation*}
\frac{d}{dt}\bar u(t)\geq(\lambda-\mu \tilde{\nu})\bar u(t)\geq-\mu\tilde{\nu}\bar u(t)\ \ \ \mbox{for all}\ t\geq T_{2},
\end{equation*}
which implies
\begin{equation*}
\bar u(t)\geq\bar u(T_{2})e^{-\mu \tilde{\nu}(t-T_{2})}\ \ \ \mbox{for all}\ t\geq T_{2}.
\end{equation*}
This along with \eqref{73} and \eqref{1144} gives
\begin{equation}\label{118}
\frac1{\bar u(t)^2}\leq C_1e^{2\mu\tilde{\nu} t}
\leq C_1e^{\frac{\rho_1}{8}t}
\ \ \ \mbox{for all}\ t\geq T_{2}.
\end{equation}
with $C_{1}>0$. Combining \eqref{72} with \eqref{118},  we can find $C_2>0$ such that
\begin{equation*}
\frac{\frac1{|\Omega|}\int_\Omega(u-\bar u)^2}{\bar u(t)^2}
\leq C_2e^{-\frac{3\rho_1}{8}t}
\ \ \ \mbox{for all}\ t\geq T_{2},
\end{equation*}
which proves \eqref{m1}.
\end{proof}

\subsection{Sharp exponential decay for $\lambda<0$}
\par\noindent\textbf{Proof of Theorem~\ref{theorem2}.}
Let $\lambda_\star$ and $T_{1}$ be as defined in Lemma \ref{l}, and let
$T_{\star, 1}:=T_{1}+1$.
Since $\lambda<0$, we have $\lambda<\lambda_\star$. The interpolation inequality in conjunction with \eqref{tt2}, \eqref{75} and  Lemma \ref{ii} provides $\beta\in(0,\alpha)$ and $C_i>0$ $(i=1,2)$ such that
\begin{equation}\label{78}
\begin{aligned}
\|\Delta v(\cdot,t)\|_{L^\infty(\Omega)}=&\|\Delta(v(\cdot,t)-\bar v(t))\|_{L^\infty(\Omega)}\\
\leq&C_1\|v(\cdot,t)-\bar v(t)\|_{C^{2+\beta}(\bar\Omega)}^{\frac{2}{2+\beta}}
\|v(\cdot,t)-\bar v(t)\|_{L^\infty(\Omega)}^{\frac{\beta}{2+\beta}}+C_1\|v(\cdot,t)-\bar v(t)\|_{L^\infty(\Omega)}\\
\leq&C_2e^{-\frac{\rho_1\beta}{12(2+\beta)}t}
\ \ \ \mbox{for all}\ t\geq T_{1}.
\end{aligned}
\end{equation}
In particular, one has
\begin{equation}\label{79}
\int_{T_{1}}^\infty\|\Delta v(\cdot,s)\|_{L^\infty(\Omega)}\,ds<\infty.
\end{equation}
By dropping a nonpositive term from the first equation in \eqref{model11}, we deduce that
\begin{equation}\label{80}
\begin{aligned}
u_t-\Delta u+\chi\nabla v\cdot\nabla u
&=(\lambda-\chi\Delta v)u-\mu u^2\\
&\leq \bigl(\lambda+\chi\|\Delta v(\cdot,t)\|_{L^\infty(\Omega)}\bigr)u
\ \ \ \mbox{in}\ \Omega\times(T_{1},\infty).
\end{aligned}
\end{equation}
Let
\begin{equation*}
H(t): =\|u(\cdot,T_{1})\|_{L^\infty(\Omega)}
\exp\left\{\lambda(t-T_{1})
+\chi\int_{T_{1}}^t\|\Delta v(\cdot,s)\|_{L^\infty(\Omega)}\,ds\right\}\ \ \ \mbox{for all}\ t\geq T_{1},
\end{equation*}
which is the solution of the associated initial-value problem
\begin{equation}\label{81}
\begin{cases}
H'(t)=\bigl(\lambda+\chi\|\Delta v(\cdot,t)\|_{L^\infty(\Omega)}\bigr)H(t),
& t>T_{1},\\[2mm]
H(T_{1})=\|u(\cdot,T_{1})\|_{L^\infty(\Omega)}.
\end{cases}
\end{equation}
Since
\begin{equation*}
u(\cdot,T_{1})\leq H(T_{1})=\|u(\cdot,T_{1})\|_{L^\infty(\Omega)},
\end{equation*}
then it follows from the parabolic comparison principle, \eqref{73} and \eqref{79} that
\begin{equation}\label{82}
\begin{aligned}
\|u(\cdot,t)\|_{L^\infty(\Omega)}\leq& H(t)\\
\leq& \|u(\cdot,T_{1})\|_{L^\infty(\Omega)}e^{-\lambda T_{1}}
\exp\!\left\{\chi\int_{T_{1}}^\infty
\|\Delta v(\cdot,s)\|_{L^\infty(\Omega)}\,ds\right\}\cdot e^{\lambda t}\\
\leq& C_{3}e^{\lambda t}
\qquad\text{for all }t\geq T_{1}
\end{aligned}
\end{equation}
with $C_{3}>0$.
As $\lambda<0$, \eqref{82} asserts that
\begin{equation}\label{86}
\int_{T_{1}}^\infty\|u(\cdot,s)\|_{L^\infty(\Omega)}\,ds<\infty.
\end{equation}

In light of \eqref{k1}, we derive
\begin{equation}\label{84}
\frac{\rm d}{{\rm d}t}\bigl(e^{-\lambda t}\bar u(t)\bigr)
=-\frac{\mu e^{-\lambda t}}{|\Omega|}\int_\Omega u^2\leq0
\ \ \ \mbox{for all}\ t\geq T_{1}.
\end{equation}
Thus the limit
\begin{equation*}
A:=\lim_{t\to\infty}e^{-\lambda t}\bar u(t)
\end{equation*}
exists and is nonnegative. We next show that $A>0$. Based on \eqref{84}, we compute
\begin{equation}\label{855}
\begin{aligned}
&\frac{\rm d}{{\rm d}t}\ln\bigl(e^{-\lambda t}\bar u(t)\bigr)
=-\mu\frac{\frac1{|\Omega|}\int_\Omega u^2}{\bar u(t)}\ \ \ \mbox{for all}\ t\geq T_{1}.
\end{aligned}
\end{equation}
An integration of \eqref{855} entails that
\begin{equation}\label{85}
\begin{aligned}
\ln\bigl(e^{-\lambda t}\bar u(t)\bigr)
=&\ln\bigl(e^{-\lambda T_{1}}\bar u(T_{1})\bigr)-\mu\int_{T_{1}}^t
\frac{\frac1{|\Omega|}\int_\Omega u^2}{\bar u(s)}ds\\
\geq&\ln\bigl(e^{-\lambda T_{1}}\bar u(T_{1})\bigr)-\mu\int_{T_{1}}^t
\|u(\cdot,s)\|_{L^\infty(\Omega)}ds\ \ \ \mbox{for all}\ t\geq T_{1},
\end{aligned}
\end{equation}
which implies
\begin{equation}\label{87}
\begin{aligned}
A\geq e^{-\lambda T_{1}}\bar u(T_{1})
\exp\left\{-\mu\int_{T_{1}}^\infty
\|u(\cdot,s)\|_{L^\infty(\Omega)}ds\right\}>0
\end{aligned}
\end{equation}
thanks to \eqref{86}. Moreover, integrating \eqref{84} in time gives
\begin{equation}\label{88}
e^{-\lambda t}\bar u(t)
=e^{-\lambda T_{1}}\bar u(T_{1})-\mu\int_{T_{1}}^t
e^{-\lambda s}\left(\frac{1}{|\Omega|}\int_\Omega u^2\right)ds\ \ \ \mbox{for all}\ t\geq T_{1},
\end{equation}
and hence
\begin{equation}\label{89}
A=e^{-\lambda T_{1}}\bar u(T_{1})-\mu\int_{T_{1}}^\infty
e^{-\lambda s}\left(\frac1{|\Omega|}\int_\Omega u^2\right)ds.
\end{equation}
Combining \eqref{88} with \eqref{89}, we infer from \eqref{82} that
\begin{equation*}
\begin{aligned}
0\leq e^{-\lambda t}\bar u(t)-A
=&\mu\int_t^\infty e^{-\lambda s}
\left(\frac1{|\Omega|}\int_\Omega u^2\right)ds\\
\leq&\mu C_3^2\int_t^\infty e^{\lambda s}\,ds
=\frac{\mu C_3^2}{-\lambda}e^{\lambda t}
\ \ \ \mbox{for all}\ t\geq T_{1}.
\end{aligned}
\end{equation*}
Consequently, there exists $C_4>0$ such that
\begin{equation}\label{97}
0\leq\bar u(t)-Ae^{\lambda t}\leq C_4e^{2\lambda t}
\ \ \ \mbox{for all}\ t\geq T_{1}.
\end{equation}

Next, we use \eqref{k1} to calculate that
\begin{equation*}
\begin{aligned}
\partial_t\bigl(e^{-\lambda t}(u-\bar u)\bigr)
=&\Delta\bigl(e^{-\lambda t}(u-\bar u)\bigr)
-\chi\nabla\cdot\bigl(e^{-\lambda t}u\nabla v\bigr)-\mu e^{-\lambda t}\left(u^2-\frac1{|\Omega|}\int_\Omega u^2\right).
\end{aligned}
\end{equation*}
Applying the variation-of-constants formula to the above identity yields
\begin{equation*}
\begin{aligned}
e^{-\lambda t}\bigl(u(\cdot,t)-\bar u(t)\bigr)
=&e^{(t-T_{1})\Delta}
\bigl[e^{-\lambda T_{1}}(u(\cdot,T_{1})-\bar u(T_{1}))\bigr]-\chi\int_{T_{1}}^t e^{(t-s)\Delta}\nabla\cdot
\bigl(e^{-\lambda s}u(\cdot,s)\nabla v(\cdot,s)\bigr)\,ds\\
&-\mu\int_{T_{1}}^t e^{(t-s)\Delta}\left[e^{-\lambda s}
\left(u(\cdot,s)^2-\frac1{|\Omega|}\int_\Omega u^2\right)\right]\,ds
\ \ \ \mbox{for all}\ t\geq T_{1},
\end{aligned}
\end{equation*}
and this leads to
\begin{equation}\label{90}
\begin{aligned}
&e^{-\lambda t}\|u(\cdot,t)-\bar u(t)\|_{L^\infty(\Omega)}\\
\leq&\underbrace{\left\|e^{(t-T_{1})\Delta}
\bigl[e^{-\lambda T_{1}}(u(\cdot,T_{1})-\bar u(T_{1}))\bigr]\right\|_{L^\infty(\Omega)}}_{:=J_3}+\underbrace{\chi\int_{T_{1}}^t \left\|e^{(t-s)\Delta}\nabla\cdot
\bigl(e^{-\lambda s}u(\cdot,s)\nabla v(\cdot,s)\bigr)\right\|_{L^\infty(\Omega)}ds}_{:=J_4}\\
&+\underbrace{\mu\int_{T_{1}}^t \left\|e^{(t-s)\Delta}\left[e^{-\lambda s}
\left(u(\cdot,s)^2-\frac1{|\Omega|}\int_\Omega u^2\right)\right]\right\|_{L^\infty(\Omega)}ds}_{:=J_5}
\ \ \ \mbox{for all}\ t\geq T_{1}.
\end{aligned}
\end{equation}
By Lemma \ref{p2}, for positive constants $C_{i}$ $(i=5, .., 10)$ we estimate that
\begin{equation}\label{91}
\begin{aligned}
J_3\leq& C_5e^{-\rho_1(t-T_{1})}
\|e^{-\lambda T_{1}}(u(\cdot,T_{1})-\bar u(T_{1}))\|_{L^\infty(\Omega)}\\
\leq&C_6e^{-\rho_1t}\ \ \ \mbox{for all}\ t\geq T_{1}
\end{aligned}
\end{equation}
and
\begin{equation}\label{92}
\begin{aligned}
J_4\leq& C_7\int_{T_{1}}^t\left(1+(t-s)^{-\frac56}\right)e^{-\rho_1(t-s)}
\|e^{-\lambda s}u(\cdot,s)\nabla v(\cdot,s)\|_{L^6(\Omega)}ds\\
\leq&C_7|\Omega|^{\frac16}\int_{T_{1}}^t\left(1+(t-s)^{-\frac56}\right)e^{-\rho_1(t-s)}\cdot e^{-\lambda s}
\|u(\cdot,s)\|_{L^\infty(\Omega)}
\|\nabla v(\cdot,s)\|_{L^\infty(\Omega)}ds\\
\leq&C_{8}\int_{T_{1}}^t\left(1+(t-s)^{-\frac56}\right)e^{-\rho_1(t-s)}\cdot e^{-\frac{\rho_1}{12}s}ds
\ \ \ \mbox{for all}\ t\geq T_{1}
\end{aligned}
\end{equation}
as well as
\begin{equation}\label{93}
\begin{aligned}
J_5\leq& C_9\int_{T_{1}}^te^{-\rho_1(t-s)}\cdot e^{-\lambda s}
\left\|u(\cdot,s)^2-\frac1{|\Omega|}\int_\Omega u^2\right\|_{L^\infty(\Omega)}ds\\
\leq&2C_9\int_{T_{1}}^te^{-\rho_1(t-s)}\cdot e^{-\lambda s}
\left\|u(\cdot,s)\right\|_{L^\infty(\Omega)}^{2}ds\\
\leq&C_{10}\int_{T_{1}}^te^{-\rho_1(t-s)}\cdot e^{\lambda s}ds
\ \ \ \mbox{for all}\ t\geq T_{1}.
\end{aligned}
\end{equation}
Here we have used \eqref{73}, \eqref{77} and \eqref{82}. Now fix any $\delta$ satisfying \eqref{ttt1}.
We combine \eqref{91}--\eqref{93} with \eqref{90} to see that
\begin{equation}\label{94}
\begin{aligned}
&e^{-\lambda t}\|u(\cdot,t)-\bar u(t)\|_{L^\infty(\Omega)}\\
\leq& C_{11}e^{-\rho_1t}+C_{11}\int_{T_{1}}^t\left(1+(t-s)^{-\frac56}\right)e^{-\rho_1(t-s)}
\cdot \left(e^{-\frac{\rho_1}{12}s}+e^{\lambda s}\right)ds\\
\leq& C_{11}e^{-\rho_1t}+2C_{11}\underbrace{\int_{T_{1}}^t\left(1+(t-s)^{-\frac56}\right)e^{-\rho_1(t-s)}
\cdot e^{-\delta s}ds}_{:=J_6}\ \ \ \mbox{for all}\ t\geq {T_{1}}.
\end{aligned}
\end{equation}
with $C_{11}>0$. For $t\geq T_{1}+1$, splitting the integral at $t-1$ gives
\begin{equation}\label{95}
\begin{aligned}
J_{6}=&e^{-\delta t}\int_{T_{1}}^t\left(1+(t-s)^{-\frac56}\right)
e^{-(\rho_1-\delta)(t-s)}\,ds\\
\leq&2e^{-\delta t}\int_{T_{1}}^{t-1}
e^{-(\rho_1-\delta)(t-s)}\,ds
+e^{-\delta t}\int_{t-1}^t\left(1+(t-s)^{-\frac56}\right)\,ds\\
\leq&\left(\frac{2}{\rho_1-\delta}+7\right)e^{-\delta t}.
\end{aligned}
\end{equation}
For $T_{1}\leq t<T_{1}+1$, $J_{6}$ can be bounded by $e^{-\delta t}\int_0^1(1+r^{-\frac56})dr=7e^{-\delta t}$. Then it follows from \eqref{94} and
\eqref{95} that
\begin{equation}\label{96}
\|u(\cdot,t)-\bar u(t)\|_{L^\infty(\Omega)}
\leq C_{12}e^{(\lambda-\delta)t}
\ \ \ \mbox{for all}\ t\geq T_{1}
\end{equation}
with $C_{12}>0$. Therefore, the combination of \eqref{97} and
\eqref{96} entails that
\begin{equation}\label{98}
\begin{aligned}
\|e^{-\lambda t}u(\cdot,t)-A\|_{L^\infty(\Omega)}\leq&e^{-\lambda t}\|u(\cdot,t)-\bar u(t)\|_{L^\infty(\Omega)}
+|e^{-\lambda t}\bar u(t)-A|\\
\leq&C_{13}\left(e^{-\delta t}+e^{\lambda t}\right)
\ \ \ \mbox{for all}\ t\geq T_{1}
\end{aligned}
\end{equation}
with some $C_{13}>0$, which together with \eqref{ttt1} establishes \eqref{ttt3} and \eqref{ttt4}.

It remains to determine the asymptotic behaviour of $v$. The fact \eqref{k3} shows that
\begin{equation}\label{100}
(v-\bar v)_t=(\Delta-1)(v-\bar v)+u-\bar u.
\end{equation}
Accordingly, the application of variation-of-constants formula and Lemma \ref{p2} gives
\begin{equation}\label{101}
\begin{aligned}
&\|v(\cdot,t)-\bar v(t)\|_{W^{1,\infty}(\Omega)}\\
\leq&\|e^{(t-T_{1})(\Delta-1)}\bigl(v(\cdot,T_{1})-\bar v(T_{1})\bigr)\|_{W^{1,\infty}(\Omega)}+\int_{T_{1}}^t \|e^{(t-s)(\Delta-1)}
\bigl(u(\cdot,s)-\bar u(s)\bigr)\|_{W^{1,\infty}(\Omega)}ds\\
\leq&C_{14}\left(1+(t-{T_{1}})^{-\frac{1}{2}}\right)e^{-(\rho_1+1)(t-T_{1})}
\bigl\|v(\cdot,T_{1})-\bar v(T_{1})\|_{L^\infty(\Omega)}\\
&+C_{15}\int_{T_1}^t
\left(1+(t-s)^{-\frac12}\right)e^{-(\rho_1+1)(t-s)}
\|u(\cdot,s)-\bar u(s)\|_{L^\infty(\Omega)}ds\\
\leq&C_{16} e^{-(\rho_1+1)t}+C_{17}\underbrace{\int_{T_1}^t
\left(1+(t-s)^{-\frac12}\right)e^{-(\rho_1+1)(t-s)}
e^{(\lambda-\delta)s}ds}_{:=J_7}\ \ \ \mbox{for all}\ t\geq T_{1}+1
\end{aligned}
\end{equation}
with $C_{i}>0$ $(i=14, .., 17)$, where we have applied \eqref{73} and \eqref{96}. We now fix $\sigma$ fulfilling \eqref{ttt2}. In view of $\delta-\lambda-\sigma>0$, we can estimate
\begin{equation*}
\begin{aligned}
J_7&=e^{-\sigma t}\int_{T_1}^t
\left(1+(t-s)^{-\frac12}\right)e^{-(\rho_1+1-\sigma)(t-s)}e^{-(\delta-\lambda-\sigma)s}\,ds
\\
&\leq2e^{-\sigma t}\int_{T_1}^{t-1}e^{-(\rho_1+1-\sigma)(t-s)}\,ds+e^{-\sigma t}\int_{t-1}^t
\left(1+(t-s)^{-\frac12}\right)ds\\
&\leq\left(\frac{2}{\rho_1+1-\sigma}+3\right)e^{-\sigma t}
\ \ \ \mbox{for all}\ t\geq T_{1}+1.
\end{aligned}
\end{equation*}
Hence \eqref{101} implies that
\begin{equation}\label{103}
\|v(\cdot,t)-\bar v(t)\|_{W^{1,\infty}(\Omega)}
\leq C_{18}e^{-\sigma t}
\ \ \ \mbox{for all}\ t\geq T_{1}+1
\end{equation}
with $C_{18}>0$.
On the other hand, solving the equation \eqref{k3} gives
\begin{equation}\label{104}
\bar v(t)=e^{-(t-T_{1})}\bar v(T_{1})
+\int_{T_{1}}^t e^{-(t-s)}\bar u(s)\,ds
\ \ \ \mbox{for all}\ t\geq T_{1}+1.
\end{equation}
We distinguish the following three cases.

If $-1<\lambda<0$, we may write
\begin{equation}\label{105}
\begin{aligned}
\bar v(t)=&e^{-(t-T_{1})}\bar v(T_{1})
+Ae^{-t}\int_{T_{1}}^t e^{(\lambda+1) s}ds+\int_{T_{1}}^t e^{-(t-s)}
\bigl(\bar u(s)-Ae^{\lambda s}\bigr)ds\\
=&\frac{A}{1+\lambda}e^{\lambda t}+e^{-(t-T_{1})}\bar v(T_{1})-\frac{A}{1+\lambda}e^{-t+(1+\lambda)T_{1}}+\int_{T_{1}}^t e^{-(t-s)}
\bigl(\bar u(s)-Ae^{\lambda s}\bigr)ds,
\end{aligned}
\end{equation}
By \eqref{ttt1} and \eqref{ttt2}, we have
$\sigma<\delta-\lambda<-2\lambda$. Hence \eqref{105} together with \eqref{73} and \eqref{97} allows us to pick $C_{19}>0$ and $C_{20}>0$ such that
\begin{equation}\label{106}
\begin{aligned}
\left|\bar v(t)-\frac{A}{1+\lambda}e^{\lambda t}\right|
&\leq C_{19}e^{-t}+C_4\int_{T_1}^t e^{-(t-s)}e^{2\lambda s}\,ds\\
&\leq C_{19}e^{-t}+C_4 e^{-\min\{1,\sigma\}t}\int_{T_1}^{\infty}e^{(2\lambda+\sigma)s}ds\\
&\leq C_{20}e^{-\min\{1,\sigma\}t}
\ \ \ \mbox{for all}\ t\geq T_{1}+1
\end{aligned}
\end{equation}
A combination of \eqref{103} and \eqref{106} now shows
\begin{equation}\label{107}
\begin{aligned}
\left\|v(\cdot,t)-\frac{A}{1+\lambda}e^{\lambda t}\right\|_{W^{1,\infty}(\Omega)}\leq&\|v(\cdot,t)-\bar v(t)\|_{W^{1,\infty}(\Omega)}+\left|\bar v(t)-\frac{A}{1+\lambda}e^{\lambda t}\right|\\
\leq&C_{18}e^{-\sigma t}
+C_{20}e^{-\min\{1,\sigma\}t}
\ \ \ \mbox{for all}\ t\geq T_{1}+1,
\end{aligned}
\end{equation}
so that \eqref{ttt5} and \eqref{ttt6} are valid thanks to $\sigma+\lambda>0$ and $-1<\lambda<0$.

If $\lambda=-1$, we rewrite \eqref{104} as
\begin{equation*}
\begin{aligned}
\bar v(t)
=&e^{-(t-T_{1})}\bar v(T_{1})+A(t-T_{1})e^{-t}
+\int_{T_{1}}^t e^{-(t-s)}\bigl(\bar u(s)-Ae^{-s}\bigr)\,ds,
\end{aligned}
\end{equation*}
so that by \eqref{97} and \eqref{73}, one has
\begin{equation}\label{108}
\begin{aligned}
|\bar v(t)-At e^{-t}|
\leq&\bigl|e^{T_{1}}\bar v(T_{1})-AT_{1}\bigr|e^{-t}
+C_4e^{-t}\int_{T_{1}}^t e^s\cdot e^{-2s}ds\\
\leq&C_{21}e^{-t}
\ \ \ \mbox{for all}\ t\geq T_{1}+1
\end{aligned}
\end{equation}
with $C_{21}>0$. Since \eqref{ttt2} ensures that
$\sigma>1$, we infer from \eqref{103} and \eqref{108} that
\begin{equation}\label{109}
\begin{aligned}
\|v(\cdot,t)-At e^{-t}\|_{W^{1,\infty}(\Omega)}\leq&\|v(\cdot,t)-\bar v(t)\|_{W^{1,\infty}(\Omega)}
+|\bar v(t)-At e^{-t}|\\
\leq&C_{18}e^{-\sigma t}+C_{21}e^{-t}\\
\leq&C_{22}e^{-t}
\ \ \ \mbox{for all}\ t\geq T_{1}+1
\end{aligned}
\end{equation}
with $C_{22}>0$. This immediately entails \eqref{ttt7} and \eqref{ttt8}.

Finally, let $\lambda<-1$. Multiplying \eqref{104} by $e^t$, we have
\begin{equation}\label{110}
e^t\bar v(t)=e^{T_{1}}\bar v(T_{1})
+\int_{T_{1}}^t e^s\bar u(s)\,ds.
\end{equation}
Due to \eqref{82} and $1+\lambda<0$, the integral in
\begin{equation}\label{111}
B:=e^{T_{1}}\bar v(T_{1})+\int_{T_{1}}^\infty e^s\bar u(s)\,ds
\end{equation}
is finite and $B>0$. Again by \eqref{82}, we combine \eqref{110} with \eqref{111} to see that
\begin{equation}\label{112}
\begin{aligned}
0\leq Be^{-t}-\bar v(t)
=&e^{-t}\int_t^\infty e^s\bar u(s)\,ds\\
\leq&C_3e^{-t}\int_t^\infty e^{(1+\lambda)s}\,ds
=\frac{C_3}{-1-\lambda}e^{\lambda t}
\ \ \ \mbox{for all}\ t\geq T_{1}.
\end{aligned}
\end{equation}
Since $\sigma>1$ also in this case, from
\eqref{103} and \eqref{112} we conclude that
\begin{equation}\label{113}
\begin{aligned}
\|v(\cdot,t)-Be^{-t} \|_{W^{1,\infty}(\Omega)}\leq&\|v(\cdot,t)-\bar v(t)\|_{W^{1,\infty}(\Omega)}
+|\bar v(t)-Be^{-t}|\\\\
\leq& C_{18}e^{-\sigma t}+\frac{C_3}{-1-\lambda}e^{\lambda t}
\ \ \ \mbox{for all}\ t\geq T_{1}+1,
\end{aligned}
\end{equation}
As $\sigma-1>0$ and $\lambda<-1$, from \eqref{113} we readily obtain \eqref{ttt9} and \eqref{ttt10}.

\subsection{Sharp algebraic decay for $\lambda=0$}
\par\noindent\textbf{Proof of Theorem~\ref{theorem3}.}
With $\lambda_{\star\star}\in(0,\lambda_{\star}]$ and $T_{2}>T_{1}$ chosen as in Lemma \ref{m}, the condition $\lambda<\lambda_{\star\star}$ holds since $\lambda=0$. Let $T_{\star, 2}:=2T_{2}+1$.  As $\bar u(t)>0$ ensured by Lemma \ref{j}, the identity \eqref{k2} allows us to calculate
\begin{equation}\label{116}
\begin{aligned}
\frac{\rm d}{{\rm d}t}\frac1{\bar u(t)}
&=-\frac{1}{\bar u^2(t)}\left\{-\mu\bar u^2(t)
-\frac{\mu}{|\Omega|}\int_\Omega(u-\bar u)^2\right\}\\
&=\mu+\mu\frac{\frac1{|\Omega|}\int_\Omega(u-\bar u)^2}
{\bar u^2(t)}
\ \ \ \mbox{for all}\ t\geq T_{2}.
\end{aligned}
\end{equation}
A direct integration of \eqref{116} leads to
\begin{equation}\label{117}
\frac1{\bar u(t)}
=\frac1{\bar u(T_{2})}+\mu(t-T_{2})
+\mu\int_{T_{2}}^t
\frac{\frac1{|\Omega|}\int_\Omega(u-\bar u(s))^2}
{\bar u^2(s)}ds
\ \ \ \mbox{for all}\ t\geq T_{2}.
\end{equation}
Recalling \eqref{m1}, we infer that the integral in
\begin{equation}\label{120}
C_\infty
:=\lim_{t\to\infty}\left(\frac1{\bar u(t)}-\mu t\right)=\frac1{\bar u(T_{2})}-\mu T_{2}
+\mu\int_{T_{2}}^\infty
\frac{\frac1{|\Omega|}\int_\Omega(u-\bar u(s))^2}
{\bar u^2(s)}ds
\end{equation}
is finite. The combination of \eqref{117} and \eqref{120} now yields
\begin{equation}\label{121}
\begin{aligned}
0\leq\mu t+C_\infty-\frac1{\bar u(t)}
&=\mu\int_t^\infty
\frac{\frac1{|\Omega|}\int_\Omega(u-\bar u(s))^2}
{\bar u^2(s)}ds\\
&\leq\mu C_2\int_t^\infty e^{-\frac{3\rho_1}{8}s}ds\\
&\leq C_3e^{-\frac{3\rho_1}{8}t}
\ \ \ \mbox{for all}\ t\geq T_{2}.
\end{aligned}
\end{equation}
with $C_{3}>0$, where we have applied \eqref{m1}. In particular, since \eqref{117} gives
\begin{equation*}
\frac{1}{\bar u(t)}\geq\frac1{\bar u(T_{2})}+\mu(t-T_{2})\ \ \ \mbox{for all}\ t\geq T_{2},
\end{equation*}
from \eqref{121} we derive that
\begin{equation}\label{122}
\begin{aligned}
\mu t+C_\infty
\geq& \frac{1}{\bar u(t)}\\
\geq& \frac{1}{(1+T_2)\bar u(T_2)}(1+T_2)+\mu(t-T_2)\\
\geq&\min\left\{\mu,\frac{1}{(1+T_2)\bar u(T_2)}\right\}(1+t)
\qquad\text{for all }t\geq T_2.
\end{aligned}
\end{equation}
Then \eqref{121} and \eqref{122} provide $C_{4}>0$ fulfilling
\begin{equation}\label{123}
\begin{aligned}
\left|\bar u(t)-\frac1{\mu t+C_\infty}\right|
&=\frac{\mu t+C_\infty-\frac1{\bar u(t)}}
{(\mu t+C_\infty)\cdot\frac1{\bar u(t)}}\\
&\leq\frac{C_4e^{-\frac{3\rho_1}{8}t}}{(1+t)^2}
\ \ \ \mbox{for all}\ t\geq T_{2}.
\end{aligned}
\end{equation}
This together with \eqref{74} entails
\begin{equation}\label{124}
\begin{aligned}
\left\|u(\cdot,t)-\frac1{\mu t+C_\infty}\right\|_{L^\infty(\Omega)}
&\leq\|u(\cdot,t)-\bar u(t)\|_{L^\infty(\Omega)}
+\left|\bar u(t)-\frac1{\mu t+C_\infty}\right|\\
&\leq C_5e^{-\frac{\rho_1}{12}t}
+\frac{C_4e^{-\frac{3\rho_1}{8}t}}{(1+t)^2}\\
&\leq C_6e^{-\frac{\rho_1}{12}t}\ \ \ \mbox{for all}\ t\geq T_{2},
\end{aligned}
\end{equation}
with $C_{5}>0$ and $C_{6}>0$, which proves \eqref{tttt1}.

Next, we turn to the second component $v$. Relying on
\begin{equation*}
\frac{\rm d}{{\rm d}t}\left(\frac1{\mu t+C_\infty}\right)
=-\frac\mu{(\mu t+C_\infty)^2}
\end{equation*}
and
\begin{equation*}
\frac{\rm d}{{\rm d}t}\left[\frac\mu{(\mu t+C_\infty)^2}\right]
=-\frac{2\mu^2}{(\mu t+C_\infty)^3},
\end{equation*}
from \eqref{k3} we infer that
\begin{equation*}
\begin{aligned}
&\frac{\rm d}{{\rm d}t}
\left[\bar v(t)-\frac1{\mu t+C_\infty}
-\frac\mu{(\mu t+C_\infty)^2}\right]
+\left[\bar v(t)-\frac1{\mu t+C_\infty}
-\frac\mu{(\mu t+C_\infty)^2}\right]\\
=&\bar u(t)-\frac1{\mu t+C_\infty}
+\frac{2\mu^2}{(\mu t+C_\infty)^3}
\ \ \ \mbox{for all}\ t\geq T_{2}.
\end{aligned}
\end{equation*}
Solving this ODI gives
\begin{equation}\label{125}
\begin{aligned}
&\bar v(t)-\frac1{\mu t+C_\infty}
-\frac\mu{(\mu t+C_\infty)^2}\\
=&\underbrace{e^{-(t-T_{2})}\left[\bar v(T_{2})-\frac1{\mu T_{2}+C_\infty}
-\frac\mu{(\mu T_{2}+C_\infty)^2}\right]}_{:=J_{8}}\\
&+\underbrace{\int_{T_{2}}^t e^{-(t-s)}
\left[\bar u(s)-\frac1{\mu s+C_\infty}
+\frac{2\mu^2}{(\mu s+C_\infty)^3}\right]ds}_{:=J_{9}}
\ \ \ \mbox{for all}\ t\geq T_{2}.
\end{aligned}
\end{equation}
Here, we may use \eqref{73} and \eqref{123} to estimate that
\begin{equation}\label{126}
\begin{aligned}
|J_{8}|\leq&C_{7}e^{-(t-T_{2})}\leq \frac{C_8}{(1+t)^3}\ \ \ \mbox{for all}\ t\geq T_{2}
\end{aligned}
\end{equation}
and
\begin{equation}\label{127}
\begin{aligned}
|J_{9}|\leq&\int_{T_{2}}^t e^{-(t-s)}
\left[\frac{C_4e^{-\frac{3\rho_1}{8}s}}{(1+s)^2}
+\frac{2\mu^2}{(\mu s+C_\infty)^3}\right]ds\\
\leq&\int_{T_{2}}^t e^{-(t-s)}
\left[\frac{C_9}{(1+s)^3}
+\frac{2\mu^2}{(\mu s+C_\infty)^3}\right]ds\\
\leq&C_{10}\int_{T_{2}}^{\frac{t}{2}}\frac{e^{-(t-s)}}{(1+s)^3}ds
+C_{10}\int_{\frac{t}{2}}^t\frac{e^{-(t-s)}}{(1+s)^3}ds\\
\leq& C_{10}e^{-\frac{t}{2}}\int_0^\infty\frac{1}{(1+s)^3}ds
+\frac{8C_{10}}{(1+t)^3}\int_{\frac{t}{2}}^t e^{-(t-s)}\,ds\\
\leq&\frac{C_{10}}{2}e^{-\frac{t}{2}}+\frac{8C_{10}}{(1+t)^3}
\leq\frac{C_{11}}{(1+t)^3}\ \ \ \mbox{for all}\ t\geq 2T_{2}+1.
\end{aligned}
\end{equation}
with $C_{i}>0$ $(i=7, ..., 11)$.
Substituting \eqref{126} and \eqref{127} into \eqref{125}, for a constant $C_{12}>0$ we have
\begin{equation*}
\left|\bar v(t)-\frac1{\mu t+C_\infty}
-\frac\mu{(\mu t+C_\infty)^2}\right|
\leq\frac{C_{12}}{(1+t)^3}
\ \ \ \mbox{for all}\ t\geq2T_{2}+1,
\end{equation*}
and hence we arrive at
\begin{equation}\label{128}
\begin{aligned}
&\left\|v(\cdot,t)-\frac1{\mu t+C_\infty}
-\frac\mu{(\mu t+C_\infty)^2}\right\|_{W^{1,\infty}(\Omega)}\\
\leq&\|v(\cdot,t)-\bar v(t)\|_{W^{1,\infty}(\Omega)}
+\left|\bar v(t)-\frac1{\mu t+C_\infty}
-\frac\mu{(\mu t+C_\infty)^2}\right|\\
\leq&C_{6}e^{-\frac{\rho_1}{12}t}+\frac{C_{12}}{(1+t)^3}
\leq\frac{C_{13}}{(1+t)^3}
\ \ \ \mbox{for all}\ t\geq2T_{2}+1
\end{aligned}
\end{equation}
thanks to \eqref{l1}. This completes the proof of \eqref{tttt2}.

It remains to prove \eqref{tttt3} and \eqref{tttt4}. We employ \eqref{124} and \eqref{128} to deduce
\begin{equation}\label{1244}
\begin{aligned}
\left\|(1+t)u(\cdot,t)-\frac1{\mu}\right\|_{L^\infty(\Omega)}
&\leq\left\|(1+t)\left(u(\cdot,t)-\frac1{\mu t+C_\infty}\right)\right\|_{L^\infty(\Omega)}
+\left|\frac{1+t}{\mu t+C_\infty}-\frac1{\mu}\right|\\
&\leq C_6(1+t)e^{-\frac{\rho_1}{12}t}+\frac{|\mu-C_\infty|}{\mu(\mu t+C_\infty)}\longrightarrow0\ \ \ \mbox{as}\ t\rightarrow\infty
\end{aligned}
\end{equation}
and
\begin{equation}\label{12444}
\begin{aligned}
&\left\|(1+t)v(\cdot,t)-\frac1\mu\right\|_{W^{1,\infty}(\Omega)}\\
\leq&
(1+t)\left\|v(\cdot,t)-\frac1{\mu t+C_\infty}
-\frac{\mu}{(\mu t+C_\infty)^2}\right\|_{W^{1,\infty}(\Omega)}
+\left|\frac{1+t}{\mu t+C_\infty}-\frac1\mu\right|
+\frac{\mu(1+t)}{(\mu t+C_\infty)^2}\\
\leq&
\frac{C_{13}}{(1+t)^2}
+\frac{|\mu-C_\infty|}{\mu(\mu t+C_\infty)}
+\frac{\mu(1+t)}{(\mu t+C_\infty)^2}
\longrightarrow 0
\qquad\text{as }t\to\infty.
\end{aligned}
\end{equation}

\subsection{Conditionally sharp exponential decay
for small $\lambda>0$}
\par\noindent\textbf{Proof of Theorem~\ref{theorem4}.}
Fix $\lambda_{\star\star}\in(0,\lambda_{\star}]$ and $T_{2}>T_{1}$ as in Lemma \ref{m}. We let $T_{\star,3}=T_{2}+1$ and $0<\lambda<\tilde{\lambda}_{0}$ with
\begin{equation}\label{keykey}
\tilde{\lambda}_{0}:=\min\left\{\lambda_{\star\star},
\frac{\rho_1}{24},\frac12\right\}.
\end{equation}
As $\bar u(t)>0$, by means of \eqref{k2} we calculate
\begin{equation*}
\begin{aligned}
\frac{\rm d}{{\rm d}t}\frac1{\bar u(t)}
&=-\frac{1}{\bar u^2(t)}\left\{\lambda\bar u(t)-\mu\bar u^2(t)
-\frac{\mu}{|\Omega|}\int_\Omega(u-\bar u)^2\right\}\\
&=-\frac{\lambda}{\bar u(t)}+\mu+\mu\frac{\frac1{|\Omega|}\int_\Omega(u-\bar u)^2}
{\bar u^2(t)}
\ \ \ \mbox{for all}\ t\geq T_{2},
\end{aligned}
\end{equation*}
which yields
\begin{equation*}
\begin{aligned}
\frac{\rm d}{{\rm d}t}
\left(\frac1{\bar u(t)}-\frac\mu\lambda\right)
+\lambda\left(\frac1{\bar u(t)}-\frac\mu\lambda\right)
&=\mu\frac{\frac1{|\Omega|}\int_\Omega(u-\bar u)^2}
{\bar u(t)^2}
\ \ \ \mbox{for all}\ t\geq T_{2}.
\end{aligned}
\end{equation*}
Solving the above ODE, we see that
\begin{equation}\label{129}
\begin{aligned}
\frac1{\bar u(t)}-\frac\mu\lambda
=e^{-\lambda t}\left\{
e^{\lambda T_{2}}\left(\frac1{\bar u(T_{2})}-\frac\mu\lambda\right)
+\mu\int_{T_{2}}^t e^{\lambda s}
\frac{\frac1{|\Omega|}\int_\Omega(u(\cdot,s)-\bar u(s))^2}
{\bar u(s)^2}\,ds\right\}\ \ \ \mbox{for all}\ t\geq T_{2}.
\end{aligned}
\end{equation}
It follows from \eqref{m1} that
\begin{equation}\label{130}
e^{\lambda t}
\frac{\frac1{|\Omega|}\int_\Omega(u-\bar u(t))^2}
{\bar u(t)^2}
\leq C_1e^{-\left(\frac{3\rho_1}{8}-\lambda\right)t}
\ \ \ \mbox{for all}\ t\geq T_{2}
\end{equation}
with $C_{1}>0$. Since $\lambda<\frac{\rho_1}{24}<\frac{3\rho_1}{8}$ by \eqref{keykey}, the right-hand side of \eqref{130} is integrable on $[T_{2},\infty)$. We may thus define
\begin{equation}\label{131}
\begin{aligned}
\widetilde{A}:=-\frac{\lambda^2}{\mu^2}\left\{
e^{\lambda T_{2}}\left(\frac1{\bar u(T_{2})}-\frac\mu\lambda\right)
+\mu\int_{T_{2}}^\infty e^{\lambda s}
\frac{\frac1{|\Omega|}\int_\Omega(u(\cdot,s)-\bar u(s))^2}
{\bar u(s)^2}\,ds\right\}.
\end{aligned}
\end{equation}
This constant is real, with no prescribed sign. Combining
\eqref{129} with \eqref{131}, we get
\begin{equation*}
\frac1{\bar u(t)}-\frac\mu\lambda
+\frac{\mu^2\widetilde{A}}{\lambda^2}e^{-\lambda t}
=-\mu e^{-\lambda t}\int_t^\infty e^{\lambda s}
\frac{\frac1{|\Omega|}\int_\Omega(u-\bar u(s))^2}
{\bar u(s)^2}\,ds\ \ \ \mbox{for all}\ t\geq T_{2}.
\end{equation*}
Therefore, we deduce from \eqref{130} that
\begin{equation}\label{132}
\begin{aligned}
\left|\frac1{\bar u(t)}-\frac\mu\lambda
+\frac{\mu^2\widetilde{A}}{\lambda^2}e^{-\lambda t}\right|
&\leq\mu C_1e^{-\lambda t}
\int_t^\infty e^{-\left(\frac{3\rho_1}{8}-\lambda\right)s}\,ds\\
&=\frac{\mu C_1}{\frac{3\rho_1}{8}-\lambda}
e^{-\frac{3\rho_1}{8}t}
\ \ \ \mbox{for all}\ t\geq T_{2},
\end{aligned}
\end{equation}
which yields the existence of $C_2>0$ satisfying
\begin{equation}\label{133}
\begin{aligned}
\left|\frac1{\bar u(t)}-\frac\mu\lambda\right|
\leq&\frac{\mu^2|\widetilde{A}|}{\lambda^2}e^{-\lambda t}
+\frac{\mu C_1}{\frac{3\rho_1}{8}-\lambda}e^{-\frac{3\rho_1}{8}t}\\
\leq& C_2e^{-\lambda t}
\ \ \ \mbox{for all}\ t\geq T_{2}.
\end{aligned}
\end{equation}
Observing the identity
\begin{equation*}
\bar u(t)-\frac\lambda\mu-\widetilde{A}e^{-\lambda t}
=-\frac{\lambda^2}{\mu^2}
\left[\frac1{\bar u(t)}-\frac\mu\lambda
+\frac{\mu^2\widetilde{A}}{\lambda^2}e^{-\lambda t}\right]
+\frac{\lambda^2}{\mu^2}\bar u(t)
\left[\frac1{\bar u(t)}-\frac\mu\lambda\right]^2,
\end{equation*}
an application of \eqref{keykey} and \eqref{132}-\eqref{133} shows that
\begin{equation}\label{134}
\begin{aligned}
\left|\bar u(t)-\frac\lambda\mu-\widetilde{A}e^{-\lambda t}\right|
\leq& C_3e^{-\frac{3\rho_1}{8}t}+C_4e^{-2\lambda t}\\
\leq& C_5e^{-2\lambda t}\ \ \ \mbox{for all}\ t\geq T_{2}
\end{aligned}
\end{equation}
with $C_{i}>0$ $(i=3, 4, 5)$. In conclusion, \eqref{134} and \eqref{l1} assert that
\begin{equation}\label{135}
\begin{aligned}
\left\|u(\cdot,t)-\frac\lambda\mu-\widetilde{A}e^{-\lambda t}
\right\|_{L^\infty(\Omega)}
&\leq\|u(\cdot,t)-\bar u(t)\|_{L^\infty(\Omega)}
+\left|\bar u(t)-\frac\lambda\mu-\widetilde{A}e^{-\lambda t}\right|\\
&\leq C_6e^{-\frac{\rho_1}{12}t}+C_5e^{-2\lambda t}\\
&\leq C_7e^{-2\lambda t}
\ \ \ \mbox{for all}\ t\geq T_{2}
\end{aligned}
\end{equation}
with $C_{6}>0$ and $C_{7}>0$, because \eqref{keykey} guarantees $2\lambda<\frac{\rho_1}{12}$. Hence, \eqref{ttttt1} and \eqref{ttttt3} become consequences of \eqref{135}

Relying on the fact
\begin{equation*}
\frac{\rm d}{{\rm d}t}
\left[\bar v(t)-\frac\lambda\mu
-\frac{\widetilde{A}}{1-\lambda}e^{-\lambda t}\right]
+\left[\bar v(t)-\frac\lambda\mu-
\frac{\widetilde{A}}{1-\lambda}e^{-\lambda t}\right]=\bar u(t)-\frac\lambda\mu-\widetilde{A}e^{-\lambda t}
\end{equation*}
due to \eqref{k3}, we may invoke the variation-of-constants formula to obtain
\begin{equation}\label{136}
\begin{aligned}
\bar v(t)-\frac\lambda\mu-\frac{\widetilde{A}}{1-\lambda}e^{-\lambda t}
=&\underbrace{e^{-(t-T_{2})}
\left[\bar v(T_{2})-\frac\lambda\mu
-\frac{\widetilde{A}}{1-\lambda}e^{-\lambda T_{2}}\right]}_{:=J_{10}}\\
&+\underbrace{\int_{T_{2}}^t e^{-(t-s)}
\left[\bar u(s)-\frac\lambda\mu-\widetilde{A}e^{-\lambda s}\right]ds}_{:=J_{11}}\ \ \ \mbox{for all}\ t\geq T_{2}.
\end{aligned}
\end{equation}
Since $2\lambda<1$  by \eqref{keykey}, we obtain from \eqref{73} and \eqref{134} that
\begin{equation}\label{137}
|J_{10}|\leq C_{8}e^{-(t-T_{2})}\ \ \ \mbox{for all}\ t\geq T_{2}
\end{equation}
and
\begin{equation}\label{138}
\begin{aligned}
|J_{11}|\leq& C_{5}e^{-t}\int_{T_{2}}^t e^{s}\cdot e^{-2\lambda s}ds\\
=&\frac{C_{5}e^{-t}\left(e^{(1-2\lambda) t}-e^{(1-2\lambda) T_{2}}\right)}
{1-2\lambda}\leq\frac{C_{5}e^{-2\lambda t}}{1-2\lambda}
\ \ \ \mbox{for all}\ t\geq T_{2}.
\end{aligned}
\end{equation}
with $C_{8}>0$, so that once more using $2\lambda<1$, the combination of \eqref{136} and
\eqref{137}-\eqref{138} implies
\begin{equation}\label{139}
\left|\bar v(t)-\frac\lambda\mu
-\frac{\widetilde{A}}{1-\lambda}e^{-\lambda t}\right|
\leq C_8e^{-(t-T_{2})}+\frac{C_5}{1-2\lambda}e^{-2\lambda t}
\leq C_9e^{-2\lambda t}
\ \ \ \mbox{for all}\ t\geq T_{2}.
\end{equation}
with $C_{9}>0$. As a consequence of \eqref{l1} and \eqref{139}, there exist $C_{10}>0$ and $C_{11}>0$ such that
\begin{equation}\label{pos:v-profile}
\begin{aligned}
&\left\|v(\cdot,t)-\frac\lambda\mu
-\frac{\widetilde{A}}{1-\lambda}e^{-\lambda t}\right\|_{W^{1,\infty}(\Omega)}\\
\leq&\|v(\cdot,t)-\bar v(t)\|_{W^{1,\infty}(\Omega)}
+\left|\bar v(t)-\frac\lambda\mu
-\frac{\widetilde{A}}{1-\lambda}e^{-\lambda t}\right|\\
\leq{}&C_{10}e^{-\frac{\rho_1}{12}t}+C_9e^{-2\lambda t}
\leq C_{11}e^{-2\lambda t}
\ \ \ \mbox{for all}\ t\geq T_{2}
\end{aligned}
\end{equation}
thanks to $2\lambda<\frac{\rho_1}{12}$, which clearly implies \eqref{ttttt2} and \eqref{ttttt4}.

\section*{Acknowledgment}
This work was supported by Shandong Provincial Natural Science Foundation (No. ZR2022JQ06),
Beijing Natural Science Foundation (Z210002) and the National Natural Science Foundation of China (12071030, 11601215, 12341102).

\end{document}